\documentclass[12pt]{article}
\usepackage{latexsym,amsfonts,amssymb}
\usepackage[dvips]{color}
\usepackage{colordvi,multicol}
\usepackage{amsmath}

\def \cal{\mathcal}

\newtheorem{thm}{Theorem}[section]

\newtheorem{lem}[thm]{Lemma}

\newtheorem{defn}[thm]{Definition}

\newtheorem{exa}[thm]{Example}

\date{}
\begin{document}
\title{\bf Periodic solutions of stochastic differential equations
driven by L\'{e}vy noises}
\author{}

\maketitle

\centerline{Xiao-Xia Guo}
\centerline{\small
School of Mathematics and Information
Sciences}
\centerline{\small Guangzhou University}
\centerline{\small  Guangzhou, 510006, China}
\centerline{\small E-mail: xxguo91@163.com}

\vskip 1cm \centerline{Wei Sun}
\centerline{\small  Department of
Mathematics and Statistics}
\centerline{\small  Concordia University}
\centerline{\small Montreal, H3G 1M8, Canada}
\centerline{\small
E-mail: wei.sun@concordia.ca}

\begin{abstract}

\noindent In this paper, we first show the well-posedness of the
SDEs  driven by L\'{e}vy noises under mild conditions. Then, we
consider the existence and uniqueness of periodic solutions of the
SDEs. To establish the ergodicity and uniqueness of periodic solutions, we
investigate the strong Feller property and the irreducibility of
the corresponding time-inhomogeneous semigroups when both small
and large jumps are allowed in the equations. Some examples are
presented to illustrate our results.
\end{abstract}




\section{Introduction}

Periodic solutions are a key concept in the theory of dynamical
systems. They
have been studied for more than a century after the pioneering work of Poincar\'e \cite{P}. Since noise is
ubiquitous in real-world systems, many people are interested in
investigating periodically varying properties of random dynamical
systems. In \cite{a8}, Khasminskii systematically studied
periodic solutions of random systems modelled by stochastic
differential equations (SDEs). But compared with the
well-developed theory for the existence of periodic solutions, the
theory of the uniqueness of periodic solutions is far from
complete.

In recent years, many works have been devoted to study the
uniqueness of periodic solutions of SDEs. Here we list some of
them which are closely related to this paper. In
\cite{a9,a11,{ra},{rb}}, Xu et al. discussed the existence and
uniqueness of periodic solutions for non-autonomous SDEs with
finite or infinite delay. In \cite{a10}, Chen et al. obtained the
existence of periodic solutions to Fokker-Planck equations through
considering the $L^2$-bounded periodic solutions in distribution
for the corresponding SDEs. In \cite{a12}, Hu and Xu presented the
existence and uniqueness theorems for periodic Markov processes on
Polish spaces. In \cite{k1}, Zhang et al.  investigated the
existence and uniqueness of periodic solutions of SDEs driven by
L\'evy processes. There are also many papers discussing
periodic solutions of stochastic biomathematical  models.
For example, Hu and Li \cite{HL} obtained the existence and uniqueness of periodic solutions of stochastic logistic equations; Zhang et al. \cite{zg1} showed  that a stochastic non-autonomous Lotka-Volterra predator-prey model with impulsive effects has a unique periodic solution, which is globally attractive.
It is worth pointing out that the above papers except for
\cite{k1}  only focused on SDEs driven by Brownian motions. Since
sudden environmental fluctuations may cause path discontinuity,
SDEs with jumps fit better the reality. We refer the reader to
\cite{bd,be,CCY, cc,cj} for some recent works on SDEs driven by
L\'{e}vy processes. The main purpose of this paper is to establish
the existence and uniqueness of periodic solutions for SDEs driven
by L\'{e}vy noises with drift, diffusion, small and
large jumps.

In the next section, we first describe the framework of this paper and
discuss the well-posedness  of global solutions to the SDEs
driven by general L\'{e}vy noises.  Then, we consider in Section 3 the
existence and uniqueness of periodic solutions.
We will follow the method of Khasminskii \cite[Theorem 3.8]{a8} to
give sufficient conditions that ensure the existence of periodic
solutions.
The main part of Section 3 is devoted to the uniqueness problem. To show the uniqueness of periodic solutions, a usual way is to
establish the asymptotic stability or global
attractivity of solutions of the SDEs under additional conditions. Completely different
from the existing methods in literature, we will investigate the strong
Feller property and the irreducibility of the time-inhomogeneous
semigroups corresponding to the SDEs. We will also extend the
uniqueness theory of invariant measures for autonomous dynamical
systems \cite[Theorem 4.2.1]{3a} to obtain the ergodicity and uniqueness
 of periodic solutions for non-autonomous SDEs driven by
L\'{e}vy noises. The conditions of our main result on periodic
solutions for SDEs with non-degenerate Gaussian noise and small/large jumps
(see Theorem \ref{lv} below) are novel and much weaker than those given in
literature.

We will use the Bismut-Elworthy-Li formula to show the strong
Feller  property of the time-inhomogeneous semigroups
corresponding to the SDEs. For the irreducibility of semigroups,
we use the method of Girsanov's transformation, which is discussed
in Da Prato and Zabczyk \cite[Theore 7.3.1]{3a}. To overcome the
difficulty caused by the jump part of L\'evy noises, we adopt the
remarkable method of Ren et. al. (see \cite{ma,cj}). For some
recent works on the uniqueness of invariant measures  for SDEs, we
recall the reader's attention to  Dong \cite{ck}, Dong and Xie
\cite{DX}, Xie \cite{ca}, Dong \cite{cc}, Xie and Zhang \cite{cj}.
In the last section of this paper, we use examples to illustrate our main
results. In particular, we will show that the stochastic
Lorenz equation and the stochastic equation of the lemniscate of Bernoulli have unique periodic solutions.

\section{Well-posedness of SDEs}
In this section, we describe the framework of the paper. We refer
the reader to  \cite{fa}  for the notation and terminology used below. Let $(\Omega, \mathcal{F},\{\mathcal{F}_t\}_{t\geq
0},\mathbb{P})$ be a complete probability space with
filtration $\{\mathcal{F}_t\}_{t\geq0}$ satisfying the usual
conditions (i.e., it is increasing, right continuous and
$\mathcal{F}_0$ contains all $\mathbb{P}$-null sets). Suppose that
$k,l, m\in \mathbb{N}$ with $k\geq m$. Denote by $\mathbb{R}_+$
the set of all non-negative real numbers. Let $\{B(t)\}_{t\ge 0}$
be a $k$-dimensional standard Brownian motion and $N$ be an
independent Poisson random measure on
$\mathbb{R}_+\times(\mathbb{R}^l-\{0\})$ with associated
compensator $\widetilde{N}$ and intensity measure $\nu$, where we
assume that $\nu$ is a L\'evy measure.

Throughout this paper, we fix a $\theta>0$. We consider the following SDE:
\begin{eqnarray}\label{a2}
dX(t)&=&b(t,X(t-))\mathrm{d}t+\sigma(t,X(t-))\mathrm{d}B(t)
+{\int_{\{|u|<1\}}}H(t,X(t-),u)\widetilde{N}
(\mathrm{d}t,\mathrm{d}u)\nonumber\\
&& +{\int_{\{|u|\geq 1\}}}G(t,X(t-),u)N
(\mathrm{d}t,\mathrm{d}u)
\end{eqnarray}
with $X(0)\in\mathcal{F}_{0}$. We assume that the coefficient
functions
$b(t,x):[0,\infty)\times\mathbb{R}^m\rightarrow\mathbb{R}^m $,
$\sigma(t,x):[0,\infty)\times\mathbb{R}^m\rightarrow
\mathbb{R}^{m\times k}$,
$H(t,x,u):[0,\infty)\times\mathbb{R}^m\times\mathbb{R}^l
\rightarrow\mathbb{R}^m$ and
$G(t,x,u):[0,\infty)\times\mathbb{R}^m\times\mathbb{R}^l
\rightarrow\mathbb{R}^m$ are all Borel measurable and satisfy
\begin{eqnarray}\label{BB}
&&b(t+\theta,x)=b(t,x),\ \ \ \ \sigma(t+\theta,x)=\sigma(t,x),\nonumber\\
&&H(t+\theta,x,u)=H(t,x,u),\ \ \ \ G(t+\theta,x,u)=G(t,x,u)
\end{eqnarray}
for any $t\ge 0$, $x\in\mathbb{R}^n$ and $u\in\mathbb{R}^l-\{0\}$.  If the large jump term is removed from
(\ref{a2}), we get the following modified SDE:
\begin{eqnarray}\label{a3}
dZ(t)=b(t,Z(t-))\mathrm{d}t+\sigma(t,Z(t-))\mathrm{d}B(t)
+{\int_{\{|u|<1\}}}H(t,Z(t-),u)\widetilde{N}
(\mathrm{d}t,\mathrm{d}u)
\end{eqnarray}
with $Z(0)=X(0)$.

We put the following assumption:

\noindent (\textbf{A1}) ~~$b(\cdot,0)~,\sigma(\cdot,0)\in
{L}^2([0,\theta);\mathbb{R}^m),
~~\int_{\{|u|<1\}}|H(\cdot,0,u)|^2\nu(\mathrm{d}u)\in
{L}^1([0,\theta);\mathbb{R}^m)$.

\noindent Hereafter we use $|x|$ to denote the Euclidean norm of a
vector $x$, use $A^T$ to denote the transpose of a matrix $A$, and
use $|A|:=\sqrt{{\rm trace}(A^TA)}$ to denote the trace norm of
$A$.

\begin{lem}\label{lem-1}
Suppose that $(\mathbf{A1})$ holds and there exists $L\in
L^1([0,\theta);\mathbb{R}_+)$ such that for any $t\in [0,\theta)$
and $x,y\in\mathbb{R}^m$,
\begin{eqnarray*}
&&| b(t,x)-b(t,y)|^2\leq L(t)|x-y|^2,\ \ \ \
|\sigma(t,x)-\sigma(t,y)|^2\leq
L(t)|x-y|^2,\nonumber\\
&&\int_{\{|u|<1\}}|H(t,x,u)-H(t,y,u)|^2\nu(\mathrm{d}u)\leq
L(t)|x-y|^2.
\end{eqnarray*}
Then, the SDE (\ref{a3}) has a unique solution $\{Z(t),t\ge
0\}$. If in addition $\mathbb{E}[|Z(0)|^2]<\infty$, then
\begin{equation}\label{gaile}
 \mathbb{E}\left[|Z(t)|^2\right]<\infty,\ \ \forall
 t\ge 0.
\end{equation}
\end{lem}

\noindent {\bf Proof.} First we assume that
$\mathbb{E}[|Z(0)|^2]<\infty$. Set $Z_0(t)=Z(0)$ for $t\ge 0$. For
$n=1,2,\dots$, define the Picard iterations
\begin{eqnarray*}
Z_n(t)&=&Z_0(t)+\int_{0}^{t}b(s,Z_{n-1}(s-))\mathrm{d}s
+\int_{0}^{t}\sigma(s,Z_{n-1}(s-))\mathrm{d}B(s)\\
&& +{\int_{0}^{t}\int_{\{|u|<1\}}}H(s,Z_{n-1}(s-),u)
\widetilde{N}(\mathrm{d}s,\mathrm{d}u). \end{eqnarray*}

Let
$$
L(t)=L(t-k\theta)\ \ \ \ {\rm for}\ t\in [k\theta,(k+1)\theta),\
k\in \mathbb{N}.
$$
By Doob's martingale inequality, for $t\ge 0$, we have
\begin{eqnarray*}
&&\mathbb{E}\left[\sup_{0\leq s\leq
t}|Z_{1}(s)-Z_{0}(s)|^2\right]\\
&\leq& 3\mathbb{E}\left[\sup_{0\leq s\leq
t}\left|\int_{0}^{s}b(v,Z_{0}(v-))
\mathrm{d}v\right|^2\right]+3\mathbb{E}\left[\sup_{0\leq s\leq
t}\left|\int_{0}^{s}
\sigma(v,Z_{0}(v-))\mathrm{d}B(v)\right|^2\right]\\
&&+ 3\mathbb{E}\left[\sup_{0\leq s\leq
t}\left|{\int_{0}^{s}\int_{\{|u|<1\}}}H(v,Z_{0}(v-),u)
\widetilde{N}(\mathrm{d}v,\mathrm{d}u)\right|^2\right]\\
&\leq&3t \mathbb{E}\left[\int_{0}^{t}|b(s,Z_{0}(s-))|^2
\mathrm{d}s\right]+12\mathbb{E}\left[\int_{0}^{t}
|\sigma(s,Z_{0}(s-))|^2\mathrm{d}s\right]\\
&&+12\mathbb{E}\left[\int_{0}^{t}\int_{\{|u|<1\}}
|H(s,Z_{0}(s-),u)|^2
\nu(\mathrm{d}u)\mathrm{d}s\right]\\
&\leq& 3t\mathbb{E}\left[\int_{0}^{t}2(|b(s,0)|^2+L(s)
|Z_{0}(s-)|^2)\mathrm{d}s\right]+12\mathbb{E}\int_{0}^{t}
2(|\sigma(s,0)|^2+L(s)|Z_{0}(s-)|^2)\mathrm{d}s\\
&&+12\mathbb{E}\int_{0}^{t}2\left(
\int_{\{|u|<1\}}|H(s,0,u)|^2\nu(\mathrm{d}u)+L(s)|Z_{0}(s-)|^2\right)
\mathrm{d}s\\
&&\leq  C_1(1+\mathbb{E}[|Z(0)|^2]),
\end{eqnarray*}
where
\begin{eqnarray*}
C_1&=&6(t+8)\int_{0}^{t}L(s)\mathrm{d}s +6t\mathbb{E}\int_{0}^{t}
|b(s,0)|^2\mathrm{d}s+24\mathbb{E}\int_{0}^{t}
|\sigma(s,0)|^2\mathrm{d}s\\
&&+24\mathbb{E}\int_{0}^{t}\int_{\{|u|<1\}}|H(s,0,u)|
^2\nu(\mathrm{d}u)ds<\infty.
\end{eqnarray*}
Similarly, we have
\begin{eqnarray*}
\mathbb{E}\left[\sup_{0\leq s\leq t}|Z_{n+1}(s)-Z_{n}(s)|^2\right]
&\leq & C_2(t)\mathbb{E}\left[\sup_{0\leq s\leq
t}|Z_{n}(s-)-Z_{n-1}(s-)|^2\right]\nonumber\\
&\leq&  (C_2(t))^n\mathbb{E}\left[\sup_{0\leq s\leq
t}|Z_{1}(s-)-Z_{0}(s-)|^2\right]\nonumber\\
&\leq& (C_2(t))^nC_1(1+\mathbb{E}[|Z(0)|^2]),
\end{eqnarray*}
where
$$C_2(t)=(3t+24)\int_{0}^{t}L(s)\mathrm{d}s.
$$

We claim that $\{Z_n(t)\}$ converges in $L^2$ for $t\geq 0$.
Indeed, for $r,n\in\mathbb{N}$ with $r<n$, we have
\begin{equation}\label{add1}
\|Z_n(t)-Z_r(t)\|_2\le\sum_{i=r+1}^{n}\|
Z_i(t)-Z_{i-1}(t)\|_2\leq\{C_1(1+\mathbb{E}[|Z(0)|^2])\}^{1/2}
\sum_{i=r+1}^{n}(C_2(t))^{i/2}.
\end{equation}
Hereafter,
$\|\cdot\|_2=\{\mathbb{E}[|\cdot|^2]\}^{1/2}$ denotes the
$L^2$-norm. We choose $\varepsilon>0$ such that
$$
(3\varepsilon+24)\int_{t}^{t+\varepsilon}L(s)\mathrm{d}s<1,\ \
\forall t\ge 0.
$$
Then, for each $0\le t\le \varepsilon$, $\{Z_n(t)\}$ is a Cauchy
sequence and hence converges to some $Z(t)\in
L^2(\Omega,\mathcal{F},\mathbb{P})$. Letting $r\rightarrow\infty$,
we obtain by (\ref{add1}) that
$$
\| Z(t)-Z_n(t)\|_2\leq\{C_1(1+\mathbb{E}[|Z(0)|^2])\}^{1/2}
\sum_{i=n+1}^{\infty}(C_2(t))^{i/2},
$$
for $n\in \mathbb{N}\cup\{0\}$. By the standard argument (cf.
\cite[Theorem 6.2.3]{fa}), we can show that $\{Z(t),0\le t\le
\varepsilon\}$ is the unique solution of the SDE (\ref{a3}) on
$[0,\varepsilon]$. Further, we  find that $\{Z(t),k\varepsilon\le
t\le (k+1)\varepsilon\}$ is the unique solution of the SDE
(\ref{a3}) on $[k\varepsilon,(k+1)\varepsilon]$ for any $k\in
\mathbb{N}$. Hence $\{Z(t),t\ge 0\}$ is the unique solution of the
SDE (\ref{a3}).

Applying the argument of  the proof of \cite[Theorem 6.2.3]{fa}, we can show
the existence and uniqueness of solutions of the SDE (\ref{a3})
for the case that $\mathbb{E}[|Z(0)|^2]=\infty$. We omit the
details here and refer the reader to \cite[Theorem
6.2.3]{fa}. \hfill\fbox

Now we put the following local Lipschitz condition.

\noindent (\textbf{A2})\ \ For each $n\in\mathbb{N}$, there exists
$L_n\in L^1([0,\theta);\mathbb{R}_+)$ such that for any $t\in
[0,\theta)$ and $x,y\in\mathbb{R}^m$ with $|x|\vee|y|\leq n$,
\begin{eqnarray*}
 &&   | b(t,x)-b(t,y)|^2\leq
  L_n(t)|x-y|^2,\ \ \ \
|\sigma(t,x)-\sigma(t,y)|^2\leq
  L_n(t)|x-y|^2,  \\
 &&   \int_{\{|u|<1\}}|H(t,x,u)-H(t,y,u)|^2\nu(\mathrm{d}u)\leq
L_n(t)|x-y|^2.
  \end{eqnarray*}

Let $C^{1,2}(\mathbb{R}_+\times\mathbb{R}^m,\mathbb{R}_+)$ be the
space of all real-valued functions $V(t,x)$ on
$\mathbb{R}_+\times\mathbb{R}^m$ which are continuously differentiable with respect to $t$ and twice continuously differentiable with respect to $x$. For $V\in C^{1,2}(\mathbb{R}_+\times\mathbb{R}^m,\mathbb{R}_+)$,
we define
\begin{eqnarray}\label{123456}
{\cal L}V(t,x)&:=&V_t(t,x)+\langle V_x(t,x), b(t,x)\rangle+
\frac{1}{2}\mathrm{ trace}(\sigma^T(t,x)V_{xx}(t,x)\sigma(t,x))\nonumber\\
&&+\int_{\{|u|<1\}}[V(t,x+H(t,x,u))-V(t,x)-
\langle V_x(t,x), H(t,x,u)\rangle]\nu(\mathrm{d}u)\nonumber\\
&&+\int_{\{|u|\ge 1\}}[V(t,x+G(t,x,u))-V(t,x)]\nu(\mathrm{d}u).
\end{eqnarray}
Hereafter we set $V_t=\frac{\partial V}{\partial t}$, $V_x=\nabla_x V=(\frac{\partial V}{\partial x_1},\dots,\frac{\partial V}{\partial x_m})$ and
$V_{xx}=(\frac{\partial^2 V}{\partial x_i\partial x_j})_{m\times m}$. A Borel measurable function $f$ on $[0,\infty)$ is said
to be locally integrable, denoted by $f\in L_{
loc}^1([0,\infty);\mathbb{R})$, if
$$
\int_0^{\tau}|f(x)|dx<\infty,\ \ \forall \tau>0.
$$
Denote by $C_0^\infty(\mathbb{R}^m)$ the space of all smooth functions on $\mathbb{R}^m$ with compact support.

Further, we make the following assumption.

\noindent ($\mathbf{H^1}$)\ \  There exist $V_1\in
C^{1,2}(\mathbb{R}_+\times\mathbb{R}^m,\mathbb{R}_+)$ and
$q_1\in L_{ loc}^1([0,\infty);\mathbb{R})$ such that
\begin{eqnarray}\label{a12}
 \lim_{|x|\rightarrow\infty}
\left[\inf_{t\in[0,\infty)}V_1(t,x)\right]=\infty,
\end{eqnarray}
and for $t\ge 0$ and $x\in\mathbb{R}^m$,
\begin{eqnarray}\label{a13}
{\cal L}V_1(t,x)\leq q_1(t).              
\end{eqnarray}

\begin{thm}\label{lem-3}
Suppose that $(\mathbf{A1})$, $(\mathbf{A2})$ and $(\mathbf{H^1})$
hold. Then, the SDE (\ref{a2}) has a unique solution $\{X(t),t\ge 0\}$.
\end{thm}
\noindent {\bf Proof.} For $n\in \mathbb{N}$, we define the
truncated functions by
\begin{eqnarray}\label{jiab}
b_n(t,x)=\left\{\begin{array}{ll}
                  b(t,x), & {\rm if}\ |x|\leq n, \\
                  b(t,\frac{n x}{|x|}), & \mathrm{elsewhere},
                \end{array}\right.
 \end{eqnarray}
 \begin{eqnarray} \label{jiabb}
\sigma_n(t,x)=\left\{\begin{array}{ll}
                  \sigma(t,x), &{\rm if}\  |x|\leq n, \\
                  \sigma(t,\frac{n x}{|x|}), & \mathrm{elsewhere},
                \end{array}\right.
\end{eqnarray}
 \begin{eqnarray}\label{jiabbb}
H_n(t,x,u)=\left\{\begin{array}{ll}
                  H(t,x,u), &{\rm if}\  |x|\leq n, \\
                  H(t,\frac{n x}{|x|},u), & \mathrm{elsewhere},
                \end{array}\right.
\end{eqnarray}
and
 \begin{eqnarray}\label{jiabbbb}
G_n(t,x,u)=\left\{\begin{array}{ll}
                  G(t,x,u), & {\rm if}\ |x|\leq n, \\
                  G(t,\frac{n x}{|x|},u), & \mathrm{elsewhere}.
                \end{array}\right.
\end{eqnarray}
Then, $b_n$, $\sigma_n$ and $H_n$ satisfy the global Lipschitz
condition and the condition (\textbf{A1}). Hence, by Lemma
\ref{lem-1}, there exists a unique solution $\{Z_n(t),t\ge 0\}$ to
the SDE
\begin{equation}\label{approx}
  \mathrm{d}Z_n(t)=b_n(t,Z_n(t-))\mathrm{d}t+\sigma_n(t,Z_n(t-))\mathrm{dB}(t)
+{\int_{\{|u|<1\}}}H_n(t,Z_n(t-),u)
\widetilde{N}(\mathrm{d}t,\mathrm{d}u)
\end{equation}
with $Z_n(0)=X(0)$.

To allow the large jump in the equation, we will use the
interlacing technique. Denote
\begin{equation}\label{BBC}
B=\{u\in\mathbb{R}^l-\{0\}:|u|<1\},\ \ \ \ B^c=\{u\in\mathbb{R}^l-\{0\}:|u|\ge1\}.
\end{equation}
Let $\{p(t)\}$ be the Poisson point process with values
in $B^c$ associated with the Poisson random measure $N(dt,du)$, i.e.,
\begin{equation}\label{BBC2}
N([0,t],A)=\#\{ p(s)\in A:s\in[0,t],A\in\mathcal{B}(B^c)\}.
\end{equation}
Then $\{p(t)\}$ is
independent of $\{Z_n(t),t\ge 0\}$, $n\in \mathbb{N}$.

Define $\tau_r:= \inf\{t>0:N([0, t]; B^c)=r\}$, which is the
$r$-th jump time of $t\mapsto N([0,t]; B^c)$. Let $\{Z(t),t\ge
0\}$ be the solution of the SDE (\ref{a3}). Define
$$X_n(t)=\left\{\begin{array}{ll}
Z_n(t), &  ~~~~\mathrm{for}~~~~ 0\le t<\tau_1,\\
Z_n(\tau_1-)+G_n(\tau_1-,Z_n(\tau_1-), p(\tau_1)), &
  ~~~~\mathrm{for}~~~~ t=\tau_1,\\
 Z_n^{(1)}(t),& ~~~~\mathrm{for}~~~~ \tau_1<t<\tau_2,\\
 Z_n^{(1)}(\tau_2-)+G_n(\tau_2-,Z_n^{(1)}(\tau_2-), p(\tau_2)),&
 ~~~~\mathrm{for}~~~~t=\tau_2,\\
 \cdots\cdots,
\end{array}\right.$$
where $\{Z_n^{(1)}(t),t\ge \tau_1\}$ is the solution of the SDE
(\ref{a3}) with $Z_n^{(1)}(\tau_1)=X_n(\tau_1)$. Then, $\{X_n(t),t\ge 0\}$ is the
unique solution of the following SDE:
\begin{eqnarray}\label{f1}
  \mathrm{d}X_n(t)&=&b_n(t,X_n(t-))\mathrm{d}t+\sigma_n(t,X_n(t-))\mathrm{dB}(t)
+{\int_{\{|u|<1\}}}H_n(t,X_n(t-),u)
\widetilde{N}(\mathrm{d}t,\mathrm{d}u)\nonumber\\
&& +{\int_{\{|u|\geq 1\}}}G_n(t,X_n(t-),u)N
(\mathrm{d}t,\mathrm{d}u)
\end{eqnarray}
with $X_n(0)=X(0)$.

For $n\in \mathbb{N}$, we define the stopping time
$$\beta_n=\inf\{t\in[0,\infty):|X_n(t)|\geq n\}.$$
For $t\in[0,\beta_n)$, we have
\begin{eqnarray*}
&&b_n(t,X_n(t))=b_{n+1}(t,X_n(t)),\ \ \ \ \sigma_n(t,X_n(t))=\sigma_{n+1}(t,X_n(t)),\\
&&H_{n}(t,X_n(t),u)=H_{n+1}(t,X_n(t),u),\ \ \ \ G_{n}(t,X_n(t),u)=G_{n+1}(t,X_n(t),u).
\end{eqnarray*}
$\{\beta_n\}$ is increasing. Hence
there exists a stopping time $\beta$ such that
$$\beta=\lim_{n\rightarrow\infty}\beta_n.
$$ Define
$$
X(t)=\lim_{n\rightarrow\infty}X_n(t),\ \ t\in [0,\beta).
$$

We now show that $\beta=\infty$ a.s.. If this is not true, then
there exist $\varepsilon>0$ and $T_1\in (0,\infty)$ such that
$$P\{\beta\leq T_1\}>2\varepsilon.$$ Hence we can find a
sufficiently large integer $n_0$ such that
\begin{eqnarray}\label{bz}
 \mathbb{P}\{\beta_n\leq T_1\}>\varepsilon,\ \ \forall n\geq n_0.
\end{eqnarray}
By It\^{o}'s formula and (\ref{a13}), we obtain that for $t\ge 0$,
\begin{eqnarray*}
&&\mathbb{E}[V_1(t\wedge\beta_n,X(t\wedge\beta_n))]\\
&=&\mathbb{E}[V_1(0,X(0))]+\mathbb{E} \left[\int_{0}^{t}
{\cal L}V_1(s,X(s\wedge\beta_n))\mathrm{d}s\right]\\
&\leq&\mathbb{E}[V_1(0,X(0))]+\int_{0}^{t}q_1(s\wedge\beta_n)\mathrm{d}s\\
&<&\infty.
\end{eqnarray*}
Thus,
$$\mathbb{E}[V_1(T_1\wedge\beta_n,X(T_1\wedge\beta_n))]<\infty,$$
which implies that
\begin{eqnarray}\label{bx}
\mathbb{E}[I_{\{\beta_n\leq T_1\}}V_1(\beta_n,X(\beta_n))]<\infty.
\end{eqnarray}

Define
$$\mu(n)=\inf\{V_1(t,x):(t,x)\in[0,\infty)\times\mathbb{R}^m,~|x|\geq n\}.$$
Then, $\lim_{n\rightarrow\infty}\mu(n)=\infty$ by the condition
(\ref{a12}). From
(\ref{bz}) and (\ref{bx}), it follows that
$$\varepsilon\mu(n)<\mu(n)\mathbb{P}\{\beta_n\leq T_1\}<\infty,$$
which results in a contradiction when $n\rightarrow\infty$.
Therefore,
\begin{equation}\label{beta}
\beta=\infty\ \ \ \ a.s.
\end{equation}
 and $\{X(t),t\ge 0\}$ is the unique
solution of the SDE (\ref{a2}) on $[0,\infty)$.\hfill\fbox

\section{Periodic solutions of SDEs driven by
L\'{e}vy noises}\label{sec2}\setcounter{equation}{0} In this
section, we will study the existence and uniqueness of periodic
solutions of the SDE (\ref{a2}). Denote by
$\mathcal{B}(\mathbb{R}^m)$ the Borel $\sigma$-algebra of
$\mathbb{R}^m$, and denote by $B_b(\mathbb{R}^m)$ (resp.
$C_b(\mathbb{R}^m)$)  the space of all real-valued bounded Borel
functions (resp. continuous and bounded functions) on
$\mathbb{R}^m$. For $f\in B_b(\mathbb{R}^m)$, we use $\|f\|_{\infty}$ to denote its supremum norm.

Recall that a stochastic process $\{X(t),t\geq 0\}$ with values in
$\mathbb{R}^m$, defined on $(\Omega,
\mathcal{F},\{\mathcal{F}_t\}_{t\geq 0},\mathbb{P})$, is called a
Markov process if, for all $A\in\mathcal{B}(\mathbb{R}^m)$ and
$0\leq s<t<\infty$,
$$\mathbb{P}\{X(t)\in A|\mathcal{F}_s\}=\mathbb{P}\{X(t)\in
A|X(s)\}.$$ We define the transition probability function of
$\{X(t),t\geq 0\}$ by
$$
P(s,x,t,A)=\mathbb{P}\{X(t)\in
A|X(s)=x\},\ \ x\in\mathbb{R}^m, A\in \mathcal{B}(\mathbb{R}^m).
$$
$\{P(s,x,t,A)\}$ defines
a semigroup of linear operators $\{P_{s,t}\}$ on $B_b(\mathbb{R}^m)$:
 $$P_{s,t}f(x):=\mathbb{E}_{s,x}[f(X(t))]:=\int_{\mathbb{R}^m}f(y)P(s,x,t,dy),\ \ x\in\mathbb{R}^m,f\in B_b(\mathbb{R}^m).$$
$\{P_{s,t}\}$ is called the Markovian transition semigroup of $\{X(t)\}$.

\begin{defn}
(i) A Markov process $\{X(t),t\geq 0\}$ is said to be
$\theta$-periodic if for any $n\in\mathbb{N}$ and any $0\le
t_1<t_2<\dots<t_n$, the joint distribution of the random variables
$X(t_1+k\theta),X(t_2+k\theta),\ldots, X(t_n+k\theta)$ is
independent of $k$ for $k\in \mathbb{N}\cup\{0\}$.
A Markovian transition semigroup $\{P_{s,t}\}$  is said to be
$\theta$-periodic if $P(s,x,t,A)=P(s+\theta,x,t+\theta,A)$ for any
$0\le s< t$, $x\in\mathbb{R}^m$ and $A\in
\mathcal{B}(\mathbb{R}^m)$. A family of probability measures $\{\mu_s, s\ge 0\}$ on $(\mathbb{R}^m,\mathcal{B}(\mathbb{R}^m))$ is said to be $\theta$-periodic with respect to
$\{P_{s,t}\}$ if
\begin{equation}\label{perio}
\mu_s(A)=\int_{\mathbb{R}^m}P(s,x,s+\theta,A)\mu_s(dx),\ \ \forall A\in \mathcal{B}(\mathbb{R}^m),\, s\ge 0.
\end{equation}

(ii) A stochastic process $\{X(t),t\geq 0\}$ with values in
$\mathbb{R}^m$ is said to be a $\theta$-periodic solution of the SDE (\ref{a2}) if it is a solution of (\ref{a2}) and is $\theta$-periodic.
\end{defn}

\begin{defn}
Let $0\le s_0<t_0<\infty$. A Markovian transition semigroup
$\{P_{s,t}\}$ is said to be regular at $(s_0,t_0)$ if all
transition probability measures $P(s_0,x,t_0,\cdot)$, $x\in\mathbb{R}^m$,
are mutually equivalent. $\{P_{s,t}\}$ is said to be Feller
(resp. strongly Feller) at $(s_0,t_0)$ if $P_{s_0,t_0}f\in
C_b(\mathbb{R}^m)$ for any $f\in C_b(\mathbb{R}^m)\
({\rm resp.}\ B_b(\mathbb{R}^m))$. $\{P_{s,t}\}$ is said to be irreducible  at
$(s_0,t_0)$ if $P(s_0,x,t_0,A)>0$ for any $x\in \mathbb{R}^m$ and
any non empty open subset $A$ of $\mathbb{R}^m$. $\{P_{s,t}\}$ is
said to be regular, Feller, strongly Feller, irreducible if it is
regular, Feller, strongly Feller, irreducible at any $(s_0,t_0)$,
respectively.
\end{defn}

\subsection{Feller and strong Feller properties of time-inhomogeneous semigroups}

Let $\{Z_n(t),t\geq 0\}$ be the solution of the SDE
(\ref{approx}). By the standard argument (cf. \cite[Theorem
6.4.5]{fa}), we can show that $\{Z_n(t),t\geq 0\}$ is a Markov
process. Further, we obtain by the
interlacing structure that the solution $\{X_n(t),t\geq 0\}$ of the SDE (\ref{f1}) is also a Markov process.
By the proof of Theorem  \ref{lem-3} and approximation, we find that
 the solution $\{X(t),t\geq
0\}$ of the SDE (\ref{a2}) is a Markov process on $\mathbb{R}^m$. In this subsection, we
will show that the transition semigroup $\{P_{s,t}\}$ of
$\{X(t),t\geq 0\}$ is Feller and strongly Feller under suitable
conditions. Let $\{X^{x}(t)\}$ be the unique solution to the SDE (\ref{a2}) with
$X^{x}(0)=x\in\mathbb{R}^m$ and let $\{X^{x}_n(t)\}$ be the unique solution to the SDE (\ref{f1}) with
$X^{x}_n(0)=x\in\mathbb{R}^m$. Denote by $\{P^n_{s,t}\}$ the
transition semigroup of $\{X_n(t),t\geq 0\}$.

We make the following assumption for the operator ${\cal L}$, which is defined in (\ref{123456}).

\noindent ($\mathbf{H_w^2}$)\ \  There exists $V_2\in C^{1,2}(\mathbb{R}_+\times\mathbb{R}^m,\mathbb{R}_+)$
such that
\begin{eqnarray}\label{jia1}
 \lim_{|x|\rightarrow\infty}
\left[\inf_{t\in[0,\infty)}V_2(t,x)\right]=\infty,
\end{eqnarray}
and
\begin{equation}\label{jia2}
\sup_{x\in\mathbb{R}^n,\,t\in[0,\infty)}{\cal L}V_2(t,x)<\infty.
\end{equation}

\begin{lem}\label{lem-4}
Suppose that $(\mathbf{A1})$, $(\mathbf{A2})$, $(\mathbf{H^1})$ and
$(\mathbf{H_w^2})$ hold.
If $\{P^n_{s,t}\}$ is  Feller for every $n\in\mathbb{N}$, then $\{P_{s,t}\}$ is Feller. If $\{P^n_{s,t}\}$ is strongly Feller for every $n\in\mathbb{N}$, then $\{P_{s,t}\}$ is strongly Feller.
\end{lem}
\noindent {\bf Proof.} To simplify notation, we only give the proof for the case that $s=0$. The proof for the case that $s>0$ is completely similar.
For $n\in \mathbb{N}$, define
$$\tau_n(x)=\inf\{t\in[0,\infty): |X_n^{x}(t)|\geq n\}.$$
By It\^{o}'s formula, we have
\begin{equation}\label{ito1}
\mathbb{E}[V_2(t\wedge\tau_n,X_n^{x}(t\wedge\tau_n))]
=\mathbb{E}[V_2(0,X_n^{x}(0))]+\mathbb{E} \left[\int_{0}^{t} {\cal
L}V_2(v,X_n^{x}(v\wedge\tau_n))\mathrm{d}v\right].
\end{equation}
By (\ref{jia1})--(\ref{ito1}), we get
\begin{eqnarray}\label{feller}
P\{\tau_n(x)<t\}&\le&\frac{V_2(0,x)+t\sup_{y\in\mathbb{R}^n,\,t\in[0,\infty)}{\cal
L}V_2(t,y)} {{\inf_{|y|>n,\,t\in[0,\infty)}}
V_2(t,y)}\nonumber\\
&\rightarrow&0\ \ \ \ {\rm as}\ n\rightarrow\infty.
\end{eqnarray}

Note that
$$X^{x}(t)=X_n^{x}(t),\ \ t<\tau_n(x).$$
Then, for $f\in B_b(\mathbb{R}^m)$ and $x,y\in\mathbb{R}^m$, we have
\begin{eqnarray}\label{Fellers}
&&|\mathbb{E}[f(X^{x}(t))]-\mathbb{E}[f(X^{y}(t))]|\nonumber\\
&=&|\mathbb{E}[f(X^{x}(t))I_{\{t<\tau_n(x)\}}]
+\mathbb{E}[f(X^{x}(t))I_{\{t\geq\tau_n(x)\}}]\nonumber\\
&& -\mathbb{E}[f(X^{y}(t))I_{\{t<\tau_n(y)\}}]
-\mathbb{E}[f(X^{y}(t))I_{\{t\geq\tau_n(y)\}}]|\nonumber\\
&=&|\mathbb{E}[f(X_n^{x}(t))]
-\mathbb{E}[f(X_n^{x}(t))I_{\{t\geq\tau_n(x)\}}]
+\mathbb{E}[f(X^{x}(t))I_{\{t\geq\tau_n(x)\}}]
-\mathbb{E}[f(X_n^{y}(t))]\nonumber\\
&&+\mathbb{E}[f(X_n^{y}(t))I_{\{t\geq\tau_n(y)\}}]
-\mathbb{E}[f(X^{y}(t))I_{\{t\geq\tau_n(y)\}}|\nonumber\\
&\leq&|\mathbb{E}[f(X_n^{x}(t))]-\mathbb{E}(f(X_n^{y}(t))]|
+|\mathbb{E}[f(X_n^{x}(t))I_{\{t\geq\tau_n(x)\}}|
+|\mathbb{E}[f(X^{x}(t))I_{\{t\geq\tau_n(x)\}}|\nonumber\\
&&+|\mathbb{E}[f(X_n^{y}(t))I_{\{t\geq\tau_n(y)\}}|+
|\mathbb{E}[f(X^{y}(t))I_{\{t\geq\tau_n(y)\}}|\nonumber\\
&\leq&|\mathbb{E}[f(X_n^{x}(t))]-\mathbb{E}[f(X_n^{y}(t))]|
+2\|{f}\|_{\infty}P(\tau_n(x)\leq t) +2\| f\|_\infty
P(\tau_n(y)\leq t).
\end{eqnarray}
Letting $n\rightarrow\infty$ and then $y\rightarrow x$ in (\ref{Fellers}), we obtain by (\ref{feller}) that $\{P_{0,t}\}$ is Feller if $\{P^n_{0,t}\}$ is Feller for every $n\in\mathbb{N}$, and $\{P_{0,t}\}$ is strongly Feller if $\{P^n_{0,t}\}$ is strongly Feller for every $n\in\mathbb{N}$.\hfill\fbox

Now we consider the Feller property of $\{P_{s,t}\}$. We need the following additional assumption.

\noindent (\textbf{B})\ \  (i) $G(t,x,u)$ is continuous in $x$ for each $t\in[0,\theta)$ and $|u|\ge 1$.

\noindent (ii) For each $n\in\mathbb{N}$,
there exist $\gamma_n>m$ and $M_n\in L^1([0,\theta);\mathbb{R}_+)$ such that for
any $t\in [0,\theta)$ and $x,y\in\mathbb{R}^m$ with
$|x|\vee|y|\leq n$,
\begin{eqnarray*}
 &&       \int_{\{|u|<1\}}|H(t,x,u)-H(t,y,u)|^{\gamma_n}\nu(\mathrm{d}u)\leq
M_n(t)|x-y|^{\gamma_n},\\
&&\int_{\{|u|<1\}}|H(t,x,u)|^{\gamma_n}\nu(\mathrm{d}u)\leq
M_n(t)(1+|x|)^{\gamma_n}.
  \end{eqnarray*}

\begin{thm}\label{thmFeller}
Suppose that $(\mathbf{A1})$, $(\mathbf{A2})$, $(\mathbf{B})$, $(\mathbf{H^1})$ and
$(\mathbf{H_w^2})$ hold. Then $\{P_{s,t}\}$ is Feller.
\end{thm}
{\bf Proof.} Under the assumption $(\mathbf{B})$, similar to the first part of the proof of \cite[Theorem 6.7.2 and Note 1 (page 402)]{fa}, we can show that $\{P^n_{s,t}\}$ is  Feller for every $n\in\mathbb{N}$. Therefore, the proof is complete by Lemma \ref{lem-4}.\hfill\fbox

\vskip 0.5cm
Next we consider the strong Feller property of $\{P_{s,t}\}$. We need the following assumptions.

\noindent (\textbf{A3})\ \ (i) $b(\cdot,0)\in
{L}^2([0,\theta);\mathbb{R}^m)$, $\sigma(\cdot,0)\in
{L}^{\infty}([0,\theta);\mathbb{R}^m)$, $\int_{\{|u|<1\}}|H(\cdot,0,u)|^2\nu(\mathrm{d}u)\in
{L}^1([0,\theta);\mathbb{R}^m)$.

\noindent (ii) For each $n\in\mathbb{N}$, there exists
$L_n\in L^{\infty}([0,\theta);\mathbb{R}_+)$ such that for any $t\in
[0,\theta)$ and $x,y\in\mathbb{R}^m$ with $|x|\vee|y|\leq n$,
\begin{eqnarray*}
 &&   | b(t,x)-b(t,y)|^2\leq
  L_n(t)|x-y|^2,\ \ \ \
|\sigma(t,x)-\sigma(t,y)|^2\leq
  L_n(t)|x-y|^2,  \\
 &&   \int_{\{|u|<1\}}|H(t,x,u)-H(t,y,u)|^2\nu(\mathrm{d}u)\leq
L_n(t)|x-y|^2.
  \end{eqnarray*}
(\textbf{A4})\ \ For any $t\in[0,\theta)$ and $x\in\mathbb{R}^m$,
$Q(t,x):=\sigma(t,x)\sigma^T(t,x)$ is invertible and \begin{equation}\label{Q}
\sup_{|x|\leq
n,\,t\in[0,\theta)}| Q^{-1}(t,x)|<\infty,\ \ \forall n\in \mathbb{N}.
\end{equation}
Obviously, (\textbf{A3}) implies (\textbf{A1}) and (\textbf{A2}).

Let $J$ be a nonnegative function in
$C_0^\infty(\mathbb{R}^m)$ satisfying
$$J(x)=0\ \ {\rm for}\ |x|\geq1\ \ \mathrm{and}\ \ \int_{\mathbb{R}^m}J(x)\mathrm{d}x=1.$$
For $\varepsilon>0$, define
$$J_\varepsilon(x)=\varepsilon^{-m}J(\varepsilon^{-1}x).$$
Let $u$ be a locally integrable function on $\mathbb{R}^m$. We define
$$
u^\varepsilon(x):=J_\varepsilon*u(x):=\int_{\mathbb{R}^m}J_\varepsilon(x-y)u(y)\mathrm{d}y.
$$

\begin{lem}\label{cor-3a}\ \  Let $\varepsilon>0$.

\noindent (i) If $u$ is a bounded function on $\mathbb{R}^m$, then
$$
\|u^\varepsilon\|_\infty\leq \|u\|_\infty.
$$

\noindent (ii) If $u$ is a continuous function on $\mathbb{R}^m$, then
$\lim_{\varepsilon\rightarrow0}u^\varepsilon(x)
=u(x)$ uniformly on  any compact subset of $\mathbb{R}^m$.

\noindent (iii) If $u$ is Lipschitz continuous on $\mathbb{R}^m$ with Lipschitz constant $L$, then
$$\| u^\varepsilon-u\|_\infty\leq L\varepsilon\ \ \ \ {\rm and}\ \ \ \ \|\nabla u^\varepsilon\|_\infty\leq L.
$$

\noindent (iv) Let $h(t,x,u):[0,\infty)\times\mathbb{R}^m\times\mathbb{R}^l
\rightarrow\mathbb{R}^m$ be Borel measurable. If there exists $\{L(t),t\ge 0\}$ such that
$$\int_{\{|u|<1\}}|h(t,x,u)-h(t,y,u)|^2
\nu(\mathrm{d}u)\leq L(t)|x-y|^2,\ \ t\ge0,\,x,y\in \mathbb{R}^m.
$$
Then,
$$\int_{\{|u|<1\}}|h^\varepsilon(t,x,u)-h^\varepsilon(t,y,u)|^2
\nu(\mathrm{d}u)\leq L(t)|x-y|^2,\ \ t\ge0,\,x,y\in \mathbb{R}^m,
$$
and
\begin{eqnarray}\label{june1}
\int_{\{|u|<1\}}|\nabla_xh^\varepsilon(t,x,u)|^2
\nu(\mathrm{d}u)\leq L(t),\ \ t\ge0,\,x\in \mathbb{R}^m.
\end{eqnarray}
\end{lem}
{\bf Proof.}
(i) The claim follows from the definition of $u^\varepsilon$.

\noindent (ii) Let $D$ be a compact subset of $\mathbb{R}^m$. We have
\begin{eqnarray}\label{JH}
|u^\varepsilon(x)-u(x)|&=&
\left|\int_{\mathbb{R}^m}
J_\varepsilon(x-y)[u(y)-u(x)]\mathrm{d}y\right|\nonumber\\
&\leq&\sup_{|y-x|\leq\varepsilon}|u(y)-u(x)|.
\end{eqnarray}
Since $u$ is uniformly continuous on $D$, (\ref{JH}) implies that
$$\lim_{\varepsilon\rightarrow0}\sup_{x\in D}|u^\varepsilon(x)-u(x)|=0.$$

\noindent (iii) This is a direct consequence of (\ref{JH}).

\noindent (iv) For $x,y\in\mathbb{R}^m$, we have
\begin{eqnarray*}
&&\int_{\{|u|<1\}}|h^\varepsilon(t,x,u)-h^\varepsilon(t,y,u)|^2\nu(du)\\
&=&\int_{\{|u|<1\}}\left|\int_{\mathbb{R}^m}(J_\varepsilon(x-v)h(t,v,u)-
J_\varepsilon(y-v)h(t,v,u))\mathrm{d}v\right|^2\nu(du)\\
&=&\int_{\{|u|<1\}}\left|\int_{\mathbb{R}^m}
(J_\varepsilon(v)h(t,x-v,u)-
J_\varepsilon(v)h(t,y-v,u))\mathrm{d}v\right|^2\nu(du)\\
&\leq&\int_{\mathbb{R}^m}
J_\varepsilon(v)\int_{\{|u|<1\}}
|h(t,x-v,u)-h(t,y-v,u)|^2\nu(du)dv\\
&\leq& L(t)|x-y|^2,
\end{eqnarray*}
which implies (\ref{june1}).
\hfill\fbox

\begin{lem}\label{cor-add}\ \  Let $\eta(t,x):[0,\infty)\times\mathbb{R}^m\rightarrow
\mathbb{R}^{m\times k}$ be Borel measurable. Suppose that

\noindent (i)  $\eta(\cdot,0)\in {L}^{\infty}([0,\theta);\mathbb{R}^m)$ and there exists $L\in L^{\infty}([0,\theta);\mathbb{R}_+)$ such that for any $x,y\in\mathbb{R}^m$,
\begin{equation}\label{approx11}
|\eta(t,x)-\eta(t,y)|^2\leq
  L(t)|x-y|^2.
\end{equation}

\noindent (ii) For any $t\in[0,\theta)$ and $x\in\mathbb{R}^m$,
$\eta(t,x)\eta^T(t,x)$ is invertible. Moreover, for any $n\in
\mathbb{N}$, there exists $\kappa_n>0$ such that
\begin{equation}\label{approx12}
\langle y,\eta(t,x)\eta^T(t,x)y\rangle\ge\kappa_n|y|^2,\ \ \ \
\forall y\in \mathbb{R}^m,\, |x|\le n,\, t\in[0,\theta).
\end{equation}
Then, for any $n\in \mathbb{N}$, if $\varepsilon$ is sufficiently small, we have
$$
\langle y,\eta^\varepsilon(t,x)(\eta^\varepsilon(t,x))^T
y\rangle\ge\frac{\kappa_n}{2}|y|^2,\ \ \ \ \forall y\in
\mathbb{R}^m,\, |x|\le n,\, t\in[0,\theta).
$$

\end{lem}
{\bf Proof.} We follow the argument of the proof of \cite[Lemma 2.2]{cc}. For $x,y\in\mathbb{R}^m$, we have
\begin{eqnarray}\label{May11}
&&\langle y,\eta^\varepsilon(t,x)(\eta^\varepsilon(t,x))^Ty\rangle\nonumber\\
&=&\sum_{i,j=1}^m\sum_{r=1}^m y_iy_j\eta_{ir}^\varepsilon(t,x)\eta_{jr}^\varepsilon(t,x)\nonumber\\
&=&\sum_{i,j=1}^m\sum_{r=1}^m y_iy_j\int_{\mathbb{R}^m}J_{\varepsilon}(x-z)\eta_{ir}(t,z)dz\int_{\mathbb{R}^m}J_{\varepsilon}(x-z')\eta_{jr}(t,z')dz'\nonumber\\
&=&\sum_{i,j=1}^m\sum_{r=1}^m y_iy_j\int_{\mathbb{R}^m}\int_{\mathbb{R}^m}J_{\varepsilon}(z)J_{\varepsilon}(z')\eta_{ir}(t,x-z)\eta_{jr}(t,x-z')dz dz'\nonumber\\
&=&\sum_{i,j=1}^m\sum_{r=1}^m y_iy_j\int_{\mathbb{R}^m}\int_{\mathbb{R}^m}J_{\varepsilon}(z)J_{\varepsilon}(z')\eta_{ir}(t,x-z)\eta_{jr}(t,x-z)dz dz'\nonumber\\
&&+\sum_{i,j=1}^m\sum_{r=1}^m y_iy_j\int_{\mathbb{R}^m}\int_{\mathbb{R}^m}J_{\varepsilon}(z)J_{\varepsilon}(z')\eta_{ir}(t,x-z)[\eta_{ir}(t,x-z')-\eta_{jr}(t,x-z)]dz dz'.\ \ \ \ \ \
\end{eqnarray}
By (\ref{approx12}), we get
\begin{eqnarray}\label{May12}
&&\sum_{i,j=1}^m\sum_{r=1}^m y_iy_j\int_{\mathbb{R}^m}\int_{\mathbb{R}^m}J_{\varepsilon}(z)J_{\varepsilon}(z')\eta_{ir}(t,x-z)\eta_{jr}(t,x-z)dz dz'\nonumber\\
&=&\int_{\mathbb{R}^m}\int_{\mathbb{R}^m}J_{\varepsilon}(z)J_{\varepsilon}(z')\langle y,\eta(t,x)\eta^T(t,x)y\rangle dz dz'\nonumber\\
&\ge& \kappa_n|y|^2,\ \ \ \ \forall y\in \mathbb{R}^m,\, |x|\le
n,\, t\in[0,\theta).
\end{eqnarray}
By (\ref{approx11}), we get
\begin{eqnarray}\label{May13}
&&\left|\int_{\mathbb{R}^m}\int_{\mathbb{R}^m}J_{\varepsilon}(z)J_{\varepsilon}(z')\eta_{ir}(t,x-z)[\eta_{ir}(t,x-z')-\eta_{jr}(t,x-z)]dz dz'\right|\nonumber\\
&\le&\int_{\mathbb{R}^m}\int_{\mathbb{R}^m}J_{\varepsilon}(z)J_{\varepsilon}(z')|\eta_{ir}(t,0)|\cdot|\eta_{ir}(t,x-z')-\eta_{jr}(t,x-z)|dz dz'\nonumber\\
&&+\int_{\mathbb{R}^m}\int_{\mathbb{R}^m}J_{\varepsilon}(z)J_{\varepsilon}(z')|\eta_{ir}(t,x-z)-\eta_{ir}(t,0)|\cdot|\eta_{ir}(t,x-z')-\eta_{jr}(t,x-z)|dz dz'\nonumber\\
&\le&\left[|\eta_{ir}(t,0)|\sqrt{L(t)}+(|x|+\varepsilon)L(t)\right]\int_{\mathbb{R}^m}\int_{\mathbb{R}^m}J_{\varepsilon}(z)J_{\varepsilon}(z')|z-z'|dz dz'\nonumber\\
&\le&2\varepsilon\left[|\eta_{ir}(t,0)|\sqrt{L(t)}+(|x|+\varepsilon)L(t)\right].
\end{eqnarray}
Therefore, the proof is complete by (\ref{May11})--(\ref{May13}).\hfill\fbox

For $n\in \mathbb{N}$, let $b_n, \sigma_n, H_n,G_n$ be defined as
in (\ref{jiab})--(\ref{jiabbbb}). Denote
$$
Q_n(t,x):=\sigma_n(t,x)\sigma_n^T(t,x).
$$
Suppose that the conditions $\mathbf{(A3)}$ and $\mathbf{(A4)}$ hold. If we replace $b, \sigma, H, Q$ with $b^\varepsilon_n,
\sigma^\varepsilon_n, H^\varepsilon_n, Q^\varepsilon_n$, then the
corresponding conditions still hold. By Lemma \ref{lem-1}, the following SDE
\begin{eqnarray}\label{zxcv}
\mathrm{d}Z^{\varepsilon}_n(t) =b^{\varepsilon}_n(t,Z^{\varepsilon}_n(t-))\mathrm{d}t
+\sigma^{\varepsilon}_n(t,Z^{\varepsilon}_n(t-))\mathrm{dB}(t)
+{\int_{\{|u|<1\}}}H^{\varepsilon}_n(t,Z^{\varepsilon}(t-),u)
\widetilde{N}(\mathrm{d}t,\mathrm{d}u)
\end{eqnarray}
 has a unique solution
$\{Z_n^{\varepsilon}(t),t\geq 0\}$ with $Z_n^{\varepsilon}(0)=x$.

\begin{lem}\label{lem-q}
Suppose that $\mathbf{(A3)}$ and $\mathbf{(A4)}$ hold. Let $T>0$.
Then, there exists a constant $M_T>0$ such that for all
$\varphi\in B_b(\mathbb{R}^m)$ and $0\le s<t\le T$,
$$
|\mathbb{E}_{s,x}[\varphi(X_n(t))]-\mathbb{E}_{s,y}[\varphi(X_n(t))]|\leq
\frac{M_T}{\sqrt{t-s}}\|\varphi\|_\infty|x-y|,\ \ \ \ \forall
x,y\in\mathbb{R}^m.
$$
\end{lem}
\noindent {\bf Proof.} We only give the proof for the case that $s=0$. The proof for the case that $s>0$ is completely similar, so we skip it.  To simplify notation, we omit the subscript ``$n$" in the
following proof.

\emph{Step 1.}
Let
$\{Z^{\varepsilon}(t),t\geq 0\}$  be the solution of the SDE
(\ref{zxcv}). Denote by $C_b^2(\mathbb{R}^m)$ the space of all continuously differentiable functions on $\mathbb{R}^m$ with bounded second-order partial derivatives. First, we prove the following Bismut-Elworthy-Li
formula: for any $\varphi\in C^2_b(\mathbb{R}^m)$, $t>0$ and
$h\in\mathbb{R}^m$, we have
\begin{eqnarray}\label{b7}
\langle
\nabla_x\mathbb{E}_{0,x}[\varphi(Z_t^{\varepsilon})],h\rangle&=&\frac{1}{t}
\mathbb{E}_{0,x}\left\{\varphi(Z^{\varepsilon}(t))
\int_{0}^t\langle
(\sigma^{\varepsilon}(t-v,Z^\varepsilon(v)))^T[\sigma^{\varepsilon}
(t-v,Z^\varepsilon(v))\right.\nonumber\\
&&\left.
\ \ \ \ \ \ \ \ \ \ \ \ \ \ \ \ (\sigma^{\varepsilon}(t-v,Z^\varepsilon(v)))^T]^{-1}\nabla_x
Z^{\varepsilon}(v)h,\mathrm{dB}(v)\rangle\right\}.
\end{eqnarray}

In fact, for $\varphi\in C^2_b(\mathbb{R}^m)$,
$\psi(t,x):=\mathbb{E}_{0,x}[\varphi(Z^{\varepsilon}(t))]$ is the
unique solution of the following equation:
$$\left\{\begin{array}{lll}
           \frac{d\psi(t,x)}{\mathrm{d}t} & = & \widetilde{L}\psi(t,x), \\
           \psi(0,x) & =& \varphi(x),\ \ \ \ t\ge 0,\,x\in\mathbb{R}^m,
         \end{array}
\right.$$
where
\begin{eqnarray*}
\widetilde{L}\psi(t,x)&:=&\langle
\psi_x(t,x),b^{\varepsilon}(t,x)\rangle+
\frac{1}{2}\mathrm{ trace}((\sigma^{\varepsilon}(t,x))^{T}\psi_{xx}(t,x)\sigma^{\varepsilon}(t,x))\\
&&+\int_{\{|u|<1\}}[\psi(t,x+H^{\varepsilon}(t,x,u))-\psi(t,x)-
\langle \psi_x(t,x),H^{\varepsilon}(t,x,u)\rangle]\nu(\mathrm{d}u).
\end{eqnarray*}
Applying It\^{o}'s formula to the process $\{\psi(t-v,
Z^{\varepsilon}(v)),v\in[0,t]\}$, we get
\begin{eqnarray}\label{b8}
\varphi(Z^{\varepsilon}(t))&=&\psi(t,x)
+{\int_{0}^{t}}\left[\frac{\partial }{\partial
v}\psi((t-v),Z^{\varepsilon}(v-))+\widetilde{L}\psi((t-v),Z^{\varepsilon}(v-))
\right]\mathrm{d}v\nonumber\\
&& +{\int_{0}^{t}}\langle
\psi_x(t-v,Z^{\varepsilon}(v-)),\sigma^{\varepsilon}(t-v,Z^{\varepsilon}(v-))
\mathrm{dB}(v)\rangle\nonumber\\
&&+{\int_{0}^{t}\int_{\{|u|<1\}}}[\psi(t-v,Z^{\varepsilon}(v-)
+H^{\varepsilon}(t-v,Z^{\varepsilon}(v-),u))
-\psi(t-v,Z^{\varepsilon}(v-))]\widetilde{N}(\mathrm{d}v,\mathrm{d}u)\nonumber\\
&=&\psi(t,x) +{\int_{0}^{t}}\langle
\psi_x(t-v,Z^{\varepsilon}(v-)),\sigma^{\varepsilon}(t-v,Z^{\varepsilon}(v-))\rangle
\mathrm{dB}(v)\nonumber\\
&&+{\int_{0}^{t}\int_{\{|u|<1\}}}[\psi(t-v,Z^{\varepsilon}(v-)
+H^{\varepsilon}(t-v,Z^{\varepsilon}(v-),u))
-\psi(t-v,Z^{\varepsilon}(v-))]\widetilde{N}(\mathrm{d}v,\mathrm{d}u).\nonumber\\
\end{eqnarray}

By Lemma \ref{cor-3a}, following the argument of \cite[Theorem 3.1]{ca}, we can show that there exists a positive constant $C_T$, which is independent of $\varepsilon$, such that
\begin{equation}\label{Bis2}
\mathbb{E}_{0,x}\{|\nabla_x Z^{\varepsilon}(v)h|^2\}\leq C_T|h|^2,\ \ \ \ \forall v\in [0,T].
\end{equation}
Note that the transition semigroup $P_{0,t}^{\varepsilon}$ of
$Z^{\varepsilon}(t)$ is given by
$P_{0,t}^{\varepsilon}\varphi(x)=\mathbb{E}_{0,x}[\varphi(Z^{\varepsilon}(t))]$.
Multiplying both sides of (\ref{b8}) by ${\int_{0}^{t}}\langle
(\sigma^{\varepsilon}(t-v,Z^\varepsilon(v)))^T(\sigma^{\varepsilon}
(t-v,Z^\varepsilon(v))(\sigma^{\varepsilon}(t-v,Z^\varepsilon(v)))^T)^{-1}\nabla_x
Z^{\varepsilon}(v)h,\mathrm{dB}(v)\rangle,$ and taking
expectation we get
\begin{eqnarray*}
&&\mathbb{E}_{0,x}\left\{\varphi(Z^{\varepsilon}(t))
{\int_{0}^{t}}\langle
(\sigma^{\varepsilon}(t-v,Z^\varepsilon(v)))^T[\sigma^{\varepsilon}
(t-v,Z^\varepsilon(v))(\sigma^{\varepsilon}(t-v,Z^\varepsilon(v)))^T]^{-1}\nabla_x Z^{\varepsilon}(v)h,\mathrm{dB}(v)\rangle\right\}\\
&=&{\int_{0}^{t}}\mathbb{E}_{0,x}\{\langle(\sigma^\varepsilon
(t-v,Z^{\varepsilon}(v)))^{T}
\psi_x(t-v,Z^{\varepsilon}(v)),(\sigma^{\varepsilon}(t-v,Z^\varepsilon(v)))^T[\sigma^{\varepsilon}
(t-v,Z^\varepsilon(v))(\sigma^{\varepsilon}(t-v,\\
&&\ \ \ \ \ \ \ \ \ \ \ \ \ \ \ \ \ \ Z^\varepsilon(v)))^T]^{-1}\nabla_x Z^{\varepsilon}(v)h\rangle\}\mathrm{d}v\\
&=&{\int_{0}^{t}}\mathbb{E}_{0,x}\{\langle
\psi_x(t-v,Z^{\varepsilon}(v)),\nabla_x Z^{\varepsilon}(v)h\rangle\}\mathrm{d}v\\
&=&{\int_{0}^{t}} \nabla_x\mathbb{E}_{0,x}\{\langle
P_{0,{t-v}}^{\varepsilon}\varphi(Z^{\varepsilon}(v)),h\rangle\}\mathrm{d}v\\
&=&t\langle \nabla_x P_{0,t}^{\varepsilon}\varphi(x),h\rangle.
 \end{eqnarray*}
Then, (\ref{b7}) holds.

By (\ref{Q}) and Lemma \ref{cor-add}, there exists $K>0$ such that
if $\varepsilon$ is sufficient small then
$$ \sup_{x\in \mathbb{R}^m,\,t\in[0,\theta)}|
[\sigma^{\varepsilon}(t,x)(\sigma^{\varepsilon}(t,x))^T]^{-1}|\le
K.
$$
Thus, we obtain by (\ref{b7}) that
\begin{eqnarray}\label{Bis1}
&&|\langle
\nabla_x\mathbb{E}_{0,x}[\varphi(Z_t^{\varepsilon})],h\rangle|^2\nonumber\\
&\leq&\frac{1}{t^2} \|\varphi\|_\infty^2
\mathbb{E}_{0,x}\left\{{\int_{0}^{t}}
|(\sigma^{\varepsilon}(t-v,Z^\varepsilon(v)))^T[\sigma^{\varepsilon}
(t-v,Z^\varepsilon(v))(\sigma^{\varepsilon}(t-v,Z^\varepsilon(v)))^T]^{-1}\nabla_x Z^{\varepsilon}(v)h|^2\mathrm{d}v\right\}\nonumber\\
&\leq&\frac{K}{t^2}\|\varphi\|_
\infty^2{\int_{0}^{t}} \mathbb{E}_{0,x}|\nabla_x
Z^{\varepsilon}(v)h|^2\mathrm{d}v.
\end{eqnarray}
Define
$$
M_T:=\sqrt{KC_T}.
$$
By (\ref{Bis2}) and (\ref{Bis1}), we get
\begin{equation}\label{h1}
|\mathbb{E}_{0,x}[\varphi(Z^{\varepsilon}(t))]-\mathbb{E}_{0,y}[\varphi(Z^{\varepsilon}(t))]|\leq
\frac{M_T}{\sqrt{t}}\|\varphi\|_\infty|x-y|,\ \ \ \ \forall
x,y\in\mathbb{R}^m.
\end{equation}

Denote by ${\rm Var}(\cdot)$ the total variation norm of a signed measure. Let $C>0$ and $t>0$. We claim that the
following conditions are equivalent:

\noindent (i) For all $\varphi\in C_b^2(\mathbb{R}^m)$ and
$x,y\in\mathbb{R}^m$,
$|P^{\varepsilon}_{0,t}\varphi(x)-P^{\varepsilon}_{0,t}\varphi(y)|\leq
C\|\varphi\|_\infty|x-y|$;

\noindent (ii) For all $\varphi\in B_b(\mathbb{R}^m)$ and
$x,y\in\mathbb{R}^m$,
$|P^{\varepsilon}_{0,t}\varphi(x)-P^{\varepsilon}_{0,t}\varphi(y)|\leq
C\|\varphi\|_\infty|x-y|$;

\noindent (iii) For all $\varphi\in B_b(\mathbb{R}^m)$ and
$x,y\in\mathbb{R}^m$,
$\mathrm{Var}(P^{\varepsilon}_{0,t}\varphi(0,x,\cdot)-P^{\varepsilon}_{0,t}\varphi(0,y,\cdot))\leq
C\|\varphi\|_\infty|x-y|$.

In fact, since each function in $B_b(\mathbb{R}^m)$ can be
approximated pointwise by functions in $C_b^2(\mathbb{R}^m)$, we have that for all $x,y\in\mathbb{R}^m$,
$$\sup_{\varphi\in\mathcal{K}_1}|P^{\varepsilon}_{0,t}\varphi(x)-P^{\varepsilon}_{0,t}\varphi(y)|
=\sup_{\varphi\in\mathcal{K}_2}|P^{\varepsilon}_{0,t}\varphi(x)-P^{\varepsilon}_{0,t}\varphi(y)|,$$
where
$$\mathcal{K}_1=\{\varphi\in
C_b(\mathbb{R}^m):\|\varphi\|_\infty<1\},
~~~~\mathcal{K}_2=\{\varphi\in
C_b^2(\mathbb{R}^m):\|\varphi\|_\infty<1\}.$$
 Hence
$$\sup_{\varphi\in\mathcal{K}_1}|P^{\varepsilon}_{0,t}\varphi(x)-P^{\varepsilon}_{0,t}\varphi(y)|=
\mathrm{Var}(P^{\varepsilon}_{0,t}\varphi(0,x,\cdot)-P^{\varepsilon}_{0,t}\varphi(0,y,\cdot)).$$
Then, (i) implies (iii). On the other hand, if (iii) holds then for all
$\varphi\in B_b(\mathbb{R}^m)$,
\begin{eqnarray*}
|P^{\varepsilon}_{0,t}\varphi(x)-P^{\varepsilon}_{0,t}\varphi(y)|&=&\left|\int_{\mathbb{R}^m}
\varphi(z)[P^{\varepsilon}(0,x,t,\mathrm{d}z)-P^{\varepsilon}(0,y,t,\mathrm{d}z)]\right|\\
&\leq&\|\varphi\|_\infty\mathrm{Var}(P^{\varepsilon}_{0,t}
\varphi(0,x,\cdot)-P^{\varepsilon}_{0,t}\varphi(0,y,\cdot))\\
&\leq&C\|\varphi\|_\infty|x-y|,
\end{eqnarray*}
which implies that (ii) holds. By the equivalence of (i)-(iii), we conclude that (\ref{h1}) holds for any $\varphi\in
B_b(\mathbb{R}^m)$.

\emph{Step 2.} In the following, we will prove that the solutions
of the SDEs (\ref{zxcv}) converge to the solution of the SDE
(\ref{approx}) in mean square, i.e.,
\begin{equation}\label{Gro}
\lim_{\varepsilon\rightarrow0}\mathbb{E}_{0,x}\left[\sup_{t\in[0,T]}|Z^\varepsilon(t)-Z(t)|^2\right]=0.
\end{equation}

In fact, by Lemma \ref{cor-3a}, we obtain that for $w\in [0,T]$,
\begin{eqnarray*}
&&\mathbb{E}_{0,x}\left[\sup_{t\in[0,w]}
|Z(t)-Z^{\varepsilon}(t)|^2\right]\\
&=&\mathbb{E}_{0,x}\left[\sup_{t\in[0,w]}
\left|\int_{0}^{t}[(b(v,Z(v-))
-b^\varepsilon(v,Z(v-)))+(b^\varepsilon(v,Z(v-))
-b^\varepsilon(v,Z^{\varepsilon}(v-)))]\mathrm{d}v\right.\right.\\
&&\left.\left.+ \int_{0}^{t}
[(\sigma(v,Z(v-))-\sigma^\varepsilon(v,Z(v-)))
+(\sigma^\varepsilon(v,Z(v-))-\sigma^\varepsilon(v,Z^{\varepsilon}(v-)))]
\mathrm{dB}(v)\right.\right.\\
&&\left.\left.+ \int_{0}^{t}\int_{\{|u|<1\}}
[(H(v,Z(v-),u)-H^\varepsilon(v,Z(v-),u))+(H^\varepsilon(v,Z(v-),u)\right.\right.\\
&&\left.\left.
-H^\varepsilon(v,Z^{\varepsilon}(v-),u))]\widetilde{N}
(\mathrm{d}v,\mathrm{d}u)\right|^2\right]\\
&\leq&6w\int_{0}^{w}L(v)\varepsilon^2dv+ 48\int_{0}^{w}L(v)\varepsilon^2dv\\
&&+(6w+48) \int_{0}^{w}
L(v)\mathbb{E}_{0,x}\left[\sup_{t\in[0,v]}
|Z(t)-Z^{\varepsilon}(t)|^2\right]\mathrm{d}v\\
&\leq&(6T+48)\varepsilon^2 \int_{0}^{T}L(v)dv +(6T+48) \int_{0}^{w}
L(v)\mathbb{E}_{0,x}\left[\sup_{t\in[0,v]}
|Z(t)-Z^{\varepsilon}(t)|^2\right]\mathrm{d}v.
\end{eqnarray*}
Define
$$
\vartheta=(6T+48)\int_{0} ^{T} L(v)dv.
$$
Then, Gronwall's inequality implies that
$$
\mathbb{E}_{0,x}\left[\sup_{t\in[0,T]} |Z(t)-Z^{\varepsilon}(t)|^2\right]\leq
\vartheta\varepsilon^2 e^{\vartheta}.
$$
Letting $\varepsilon\rightarrow0$, we obtain (\ref{Gro}).

Finally, combining (\ref{h1}) and (\ref{Gro}), we obtain that for any $\varphi\in B_b(\mathbb{R}^m)$, $t\in (0,T]$ and $x,y\in\mathbb{R}^m$,
\begin{equation}\label{nan}
\left|\mathbb{E}_{0,x}[\varphi(Z(t))]-\mathbb{E}_{0,y}[\varphi(Z(t))]\right|\leq
\frac{M_T}{\sqrt{t}}\|\varphi\|_\infty|x-y|.
\end{equation}

\emph{Step 3.} Let $B, B^c$ and the  Poisson point process $\{p(t)\}$ be defined as in (\ref{BBC}) and (\ref{BBC2}), respectively.
Define $\tau_1:= \inf\{t>0:N([0, t]; B^c)=1\}$, which is the
first jump time of $t\mapsto N([0,t]; B^c)$. Let $\{Z(t),t\ge
0\}$ be the solution of the SDE (\ref{approx}). Then $\{p(t)\}$ is
independent of $\{Z(t),t\ge 0\}$.

Denote by $P^Z_{0,t}$ and $P^X_{0,t}$ the transition semigroups of
$Z(t)$ and $X(t)$, respectively. Then, for $\varphi\in B_b(\mathbb{R}^m)$, $t\ge 0$ and $x\in\mathbb{R}^m$, we have
$$
P^Z_{0,t}\varphi(x)=\mathbb{E}_{0,x}[\varphi(Z(t))],\ \ \ \ P^X_{0,t}\varphi(x)=\mathbb{E}_{0,x}[\varphi(X(t))].
$$
Note that
\begin{eqnarray}\label{nan2}
P^X_{0,t}\varphi(x)&=&\mathbb{E}_{0,x}[\varphi(X(t));t<\tau_1]+
\mathbb{E}_{0,x}[\varphi(X(t));t\ge\tau_1]\nonumber\\
&=&\mathbb{E}_{0,x}[\varphi(Z(t));t<\tau_1]+\mathbb{E}_{0,x}\{I_{\{\tau_1\le t\}}\mathbb{E}_{\tau_1,X(\tau_1)}
[\varphi(X(t))]\}\nonumber\\
&=&e^{-\nu(B^c)t}P^Z_{0,t}\varphi(x)\nonumber\\
&&+
\int_{0}^t\int_{\{|u|\ge 1\}}\int_{\mathbb{R}^m}
e^{-\nu(B^c)v}P^X_{v,t}\varphi(w+G(v,w,u))P^Z_{0,v}(x,dw)\nu(du)dv.
\end{eqnarray}
By (\ref{h1}) and (\ref{nan2}), we obtain that for $t\in (0,T]$,
\begin{eqnarray*}
&&\left|\mathbb{E}_{0,x}[\varphi(X(t))]-\mathbb{E}_{0,y}[\varphi(X(t))]\right|\\
&\leq&\left|\mathbb{E}_{0,x}[\varphi(Z(t))]-\mathbb{E}_{0,y}[\varphi(Z(t))]\right|\\
&&+\left|\int_{0}^t\int_{\{|u|\ge 1\}}\int_{\mathbb{R}^m}
e^{-\nu(B^c)v}P^X_{v,t}\varphi(w+G(v,w,u))[P^Z_{0,v}(x,dw)-P^Z_{0,v}(y,dw)]\nu(du)dv\right|\\
&\leq&
\frac{M_T}{\sqrt{t}}\|\varphi\|_\infty|x-y|+\int_{0}^t\nu(B^c)e^{-\nu(B^c)v}\frac{M_T}{\sqrt{v}}\|\varphi\|_\infty|x-y|dv\\
&=&M_T\|\varphi\|_\infty|x-y|\left(\frac{1}{\sqrt{t}}+2\nu(B^c)\sqrt{t}\right)\\
&\le&\frac{M_T+2\nu(B^c)T}{\sqrt{t}}\|\varphi\|_\infty|x-y|.
\end{eqnarray*}
Therefore, the proof is complete.\hfill\fbox

\begin{thm}\label{thm-s}
Suppose that $\mathbf{(A3)}$, $\mathbf{(A4)}$, $(\mathbf{H^1})$ and
$(\mathbf{H_w^2})$ hold. Then $\{P_{s,t}\}$ is strongly Feller.
\end{thm}
{\bf Proof.} By Lemma \ref{lem-q}, we know that $\{P^n_{s,t}\}$ is strongly Feller for every $n\in\mathbb{N}$. Hence Lemma \ref{lem-4} implies that $\{P_{s,t}\}$ is strongly Feller.\hfill\fbox

\subsection{Irreducibility of time-inhomogeneous semigroups}

Let $\{X(t),t\geq
0\}$ be the solution of the SDE (\ref{a2}). In this subsection, we
will show that the transition semigroup $\{P_{s,t}\}$ of
$\{X(t),t\geq 0\}$ is irreducible. Denote by $B_{b,loc}(\mathbb{R}_+)$ the set of all locally bounded Borel measurable functions on $\mathbb{R}_+$.
We make the following assumption.

\noindent $(\mathbf{H^3})$\ \  (i) There exists $W(t,x):[0,\infty)\times\mathbb{R}^m\rightarrow\mathbb{R}^m$ which is locally bounded and for each $n\in\mathbb{N}$, there exists
$R_n\in L_{loc}^1([0,\infty);\mathbb{R}_+)$ such that for any $t\in
[0,\infty)$ and $x,y\in\mathbb{R}^m$ with $|x|\vee|y|\leq n$,
$$
   |W(t,x)-W(t,y)|^2\leq
  R_n(t)|x-y|^2.
  $$

  \noindent (ii) There exist $q_3,q\in B_{b,loc}(\mathbb{R}_+)$ and
$U(t,x):\mathbb{R}_+\times \mathbb{R}_+\rightarrow(0,\infty)$ which is increasing with respect to $x$.

  \noindent (iii) There exists $V_3\in
C^{1,2}([0,\infty)\times\mathbb{R}^m,\mathbb{R}_+)$  such that
$$ \lim_{|x|\rightarrow\infty}
\left[\inf_{t\in[0,\infty)}V_3(t,x)\right]=\infty,
$$
and for $t\ge 0$ and $x\in\mathbb{R}^m$,
\begin{eqnarray}\label{mya}
{\cal L}V_3(t,x)\leq q_3(t),
\end{eqnarray}
\begin{eqnarray}\label{W}
U(t,|x|)\le V_3(t,x)\le\langle W,\nabla_x V_3(t,x)\rangle+q(t).
\end{eqnarray}

\begin{thm}\label{thm-6} Suppose that $(\mathbf{A1})$, $(\mathbf{A2})$, $(\mathbf{A4})$ and $(\mathbf{H^3})$ hold.
Then $\{P_{s,t}\}$ is irreducible.
\end{thm}
{\bf Proof.}  To simplify notation, we only give the proof for the case that $s=0$. The proof for the case that $s>0$ is completely similar.
 Let $x,y\in\mathbb{R}^m$ with $x\not=y$ and $T>0$. For $r\in\mathbb{N}$, we consider the
following SDE:
\begin{eqnarray}\label{zb}
  dX^r(t) & = & \left[b(t,X^r(t-))-rW(t,X^r(t-)-y)\right]dt
    +\sigma(t,X^r(t-))dB(t)\nonumber\\
    &&+{\int_{\{|u|<1\}}}H(t,X^r(t-),u)\widetilde{N}
(\mathrm{d}t,\mathrm{d}u)+{\int_{\{|u|\ge1\}}}G(t,X^r(t-),u){N}
(\mathrm{d}t,\mathrm{d}u).\ \ \ \
\end{eqnarray}
Note that the generator ${\cal L}^W$ of the SDE ({\ref{zb}) is given by
$$
{\cal L}^Wf(t,x)={\cal L}f(t,x)-r\langle W,\nabla_x f\rangle(t,x).
$$
By the conditions $(\mathbf{A1})$, $(\mathbf{A2})$ and $(\mathbf{H^3})$, following the argument of the proof of Theorem \ref{lem-3}, we can show that the SDE ({\ref{zb}) has a unique solution $\{X^{r}(t),t\ge 0\}$ with $X^r(0)=x$.

By It\^{o}'s formula and (\ref{W}), we get
\begin{eqnarray*}
  &&\mathbb{E}[e^{rt}V_3(t,X^r(t)-y)]\\
    &=&  \mathbb{E}[V_3(0,x-y)]+\mathbb{E}\left[\int_{0}^{t}
  e^{rv}{\cal L}^WV_3(v,X^r(v-)-y)\mathrm{d}v\right]
     \\
  &&+\mathbb{E}
  \left[\int_{0}^{t}re^{rv}V_3(s,X^r(v-)-y)\mathrm{d}v\right]\\
    &\leq&  \mathbb{E}[V_3(0,x-y)]+\mathbb{E}\left[\int_{0}^{t}
e^{rv}
[-rV_3(v,X^r(v-)-y)+rq(v)+q_3(v)]\mathrm{d}v\right] \\
&&+\mathbb{E}\left[\int_{0}^{t}r
e^{rv}V_3(v,X^r(v-)-y)
\mathrm{d}v\right]\\
&=&  \mathbb{E}[V_3(0,x-y))]+\int_{0}^{t}
(rq(v)+q_3(v))e^{rv}\mathrm{d}v.
\end{eqnarray*}
Hence
$$
  \mathbb{E}[V_3(t,X^r(t)-y)]\leq
  \frac{\mathbb{E}[V_3(0,x-y)]}{e^{rt}}+\frac{(\sup_{0\leq
v\le t} |rq(v)+q_3(v)|)(1-e^{-rt})}{r}.
$$
Then, for any
$0<a<|x-y|$, there exists $r_a\in \mathbb{N}$ such that
\begin{eqnarray}\label{zv}
&&\mathbb{P}(|X^{r_a}(T)-y|\ge a)\nonumber\\
&\leq&
\mathbb{P}\{U(T,|X^{r_a}(T)-y|)\geq U(T,a)\} \nonumber \\
 &\leq&  \mathbb{P}\{V_3(T,X^{r_a}(T)-y)\geq U(T,a)\}\nonumber \nonumber \\
 &\leq& \frac{\mathbb{E}[V_3(T,X^{r_a}(T)-y)]}{U(T,a)}\nonumber \\
 &\leq& \frac{\mathbb{E}[V_3(0,x-y)]}{e^{r_aT}U(T,a)}
 +\frac{(\sup_{0\leq
v\le T} |r_aq(v)+q_3(v)|)(1-e^{-r_aT})}{r_a U(T,a)}\nonumber \\
&<&\frac{1}{2}.
\end{eqnarray}

For $K\in \mathbb{N}$, define
$$\tau_K:=\inf\{t:|X^{r_a}(t)|\geq K\}.$$
Similar to (\ref{beta}), we can show that there exists $K\in \mathbb{N}$ such that
$$
\mathbb{P}(\tau_K\leq T)<\frac{1}{2},
$$
which together with (\ref{zv}) implies that
\begin{equation}\label{vt}
\mathbb{P}(\tau_K\leq T)+\mathbb{P}(|X^{r_a}(T)-y|\ge a)<1.
\end{equation}
Define
$$\alpha(t):=-r_a(\sigma(t,X^{r_a}(t)))^T[\sigma(t,X^{r_a}(t))(\sigma(t,X^{r_a}(t)))^T]^{-1}W(t,X^{r_a}(t)-y),
$$
$$
\widetilde{B}(t):=B(t)+\int_{0}^{t\wedge \tau_K}\alpha(v)dv,$$
and
$$\mathfrak{M}(t)=\mathrm{exp}
\left(\int_{0}^{t\wedge\tau_K}\alpha(v)dB(v)
-\frac{1}{2}\int_{0}^{t\wedge\tau_K}|\alpha(v)|^2dv\right).$$
By $(\mathbf{A1})$, $(\mathbf{A2})$, $(\mathbf{A4})$ and $(\mathbf{H^3})$, we get
$$
\mathbb{E}\left\{\mathrm{exp}
\left(\frac{1}{2}\int_{0}^{t\wedge\tau_K}|\alpha(v)|^2dv\right)\right\}<\infty,
$$
i.e., Novikov's condition is satisfied. Then $\{\mathfrak{M}(t)\}$ is a martingale. Thus,  by Girsanov's theorem, under
the new probability measure
$\mathbb{Q}=\mathfrak{M}(t)\mathbb{P}$, $\widetilde{B}(t)$ is
still a Brownian motion, and $N(dt,du)$ is a Poisson random
measure with the same compensator $\nu(du)dt$.

By (\ref{vt}), we have
\begin{equation}\label{QQ}
\mathbb{Q}(\{\tau_K\leq T\}\cup\{|X^{r_a}(T)-y|\ge a\})<1.
\end{equation}
Note that $X^{r_a}(t)$ also solves the following SDE:
\begin{eqnarray*}
 X(t\wedge\tau_K) &=& x+\int_{0}^{t\wedge\tau_K}b(v,X(v-))ds
+\int_{0}^{t\wedge\tau_K}\sigma(v,X(v-))d\widetilde{B}(v) \nonumber\\
   &
   &+\int_{0}^{t\wedge\tau_K}\int_{\{|u|<1\}}H(v,X(v-),u)\widetilde{N}(dv,du)\nonumber\\
   &
   &+\int_{0}^{t\wedge\tau_K}\int_{\{|u|\ge1\}}G(v,X(v-),u){N}(dv,du).
   \end{eqnarray*}
  Set $\eta_K:=\inf\{t:|X(t)|\geq K\}$. By the weak uniqueness of the solutions of the SDE
(\ref{a2}), we know that the law of
$\{(X(t)I_{\{t<\eta_K\}}\}$ under $\mathbb{P}$ is the
same as that of $\{(X^{r_a}(t)I_{\{t<\tau_K\}}\}$
under $\mathbb{Q}$. Hence we obtain by (\ref{QQ}) that
\begin{eqnarray*}
    \mathbb{P}(|X(T)-y|\ge a) & \leq & \mathbb{P}(\{\eta_K\leq T\}
    \cup\{\eta_K> T,|X(T)-y|\ge a\}) \\
     & = & \mathbb{Q}(\{\tau_K\leq T\}
    \cup\{\tau_K> T,|X^{r_a}(T)-y|\ge a\}) \\
     & \leq & \mathbb{Q}(\{\tau_K\leq T\}
    \cup\{|X^{r_a}(T)-y|\ge a\})\\
    & <&1,
  \end{eqnarray*}
which implies that
$$
  \mathbb{P}(|X(T)-y|< a)>0.
$$
Since $a,x,y, T$ are arbitrary, we conclude that $\{P_{s,t}\}$ is irreducible.\hfill\fbox

\subsection{Existence of periodic solutions}
\begin{lem}\label{lem-3a}
Let $\{X(t),t\geq 0\}$ be the unique solution of the SDE
(\ref{a2}). Then its transition semigroup $\{P_{s,t}\}$ is $\theta$-periodic.
\end{lem}
\noindent {\bf Proof.} Define
$$\overline{B}(t)=B(t+\theta)-B(\theta).
$$
Then, we obtain by (\ref{a2}) and (\ref{BB}) that
\begin{eqnarray*}
&&X(t+\theta)\\
&=& X(0)+{\int_{0}^{\theta}}b(s,X(s-))ds
+{\int_{0}^{\theta}}\sigma(s,X(s-))\mathrm{dB}(s)\\
&&+ \int_{0}^{\theta}\int_{\{|u|<1\}}H(s,X(s-),u)
\widetilde{N}(\mathrm{d}s,\mathrm{d}u)
 +\int_{0}^{\theta}\int_{\{|u|\geq
1\}}G(s,X(s-),u)N(\mathrm{d}s,\mathrm{d}u)\\
&&+\int_{\theta}^{t+\theta}b(s,X(s-))ds
+\int_{\theta}^{t+\theta}\sigma(s,X(s-))\mathrm{dB}(s)\\
&&+ \int_{\theta}^{t+\theta}\int_{\{|u|<1\}}H(s,X(s-),u)
\widetilde{N}(\mathrm{d}s,\mathrm{d}u)+\int_{\theta}^{t+\theta}\int_{\{|u|\geq
1\}}G(s,X(s-),u)
N(\mathrm{d}s,\mathrm{d}u)\\
&=& X(\theta)+\int_{\theta}^{t+\theta}b(s,X(s-))ds
+\int_{\theta}^{t+\theta}\sigma(s,X(s-))\mathrm{dB}(s)\\
&&+ \int_{\theta}^{t+\theta}\int_{\{|u|<1\}}H(s,X(s-),u)
\widetilde{N}(\mathrm{d}s,\mathrm{d}u)
+\int_{\theta}^{t+\theta}\int_{\{|u|\geq 1\}}G(s,X(s-),u)
N(\mathrm{d}s,\mathrm{d}u)\\
&=& X(\theta)+\int_{0}^{t}b(r+\theta,X(\theta +r-))dr
+\int_{0}^{t}\sigma(\theta +r,X(\theta+ r-))\mathrm{d\overline{B}}(r)\\
&&+ \int_{0}^{t}\int_{\{|u|<1\}}H(\theta +r,X(\theta+r-),u)
\widetilde{N}(\mathrm{d}r,\mathrm{d}u)\\
&&+\int_{0}^{t}\int_{\{|u|\geq 1\}}G(\theta +r,X(\theta +r-),u)
N(\mathrm{d}r,\mathrm{d}u)\\
&=& X(\theta)+\int_{0}^{t}b(r,X(\theta +r-))dr
+\int_{0}^{t}\sigma(r,X(\theta +r-))\mathrm{d\overline{B}}(r)\\
&&+ \int_{0}^{t}\int_{\{|u|<1\}}H(r,X(\theta+ r-),u)
\widetilde{N}(\mathrm{d}r,\mathrm{d}u)
+\int_{0}^{t}\int_{\{|u|\geq 1\}}G(r,X(\theta +r-),u)
N(\mathrm{d}r,\mathrm{d}u).
\end{eqnarray*}
Hence
\begin{eqnarray}\label{app2}
\mathrm{d}X(t+\theta)&=&b(t,X(t+\theta-))\mathrm{d}t
+\sigma(t,X(t+\theta-))\mathrm{d\overline{B}}(t)
+\int_{\{|u|<1\}}H(r,X(r+\theta-),u)
\widetilde{N}(\mathrm{d}r,\mathrm{d}u)\nonumber\\
&& +\int_{\{|u|\geq 1\}}G(r,X(r+\theta-),u)
N(\mathrm{d}r,\mathrm{d}u).
\end{eqnarray}

By (\ref{app2}), we find that  $\{X(t+\theta),t\geq 0\}$ is  a
weak solution of  the SDE (\ref{a2}). From the weak uniqueness of solutions, we
know that $\{X(t),t\geq 0\}$ and $\{X(t+\theta),t\geq 0\}$ have
the same distribution. Therefore,
$$P(s,x,t,A)=P(s+\theta,x,t+\theta,A),\ \ \ \ \forall 0\le s< t,x\in\mathbb{R}^m,A\in \mathcal{B}(\mathbb{R}^m).
$$
\hfill\fbox

We make the following assumption for the operator ${\cal L}$, which is defined in (\ref{123456}).

\noindent ($\mathbf{H^2}$)\ \  There exists $V_2\in C^{1,2}(\mathbb{R}_+\times\mathbb{R}^m,\mathbb{R}_+)$
such that
\begin{eqnarray}\label{bq1}
 \lim_{|x|\rightarrow\infty}
\left[\inf_{t\in[0,\infty)}V_2(t,x)\right]=\infty,
\end{eqnarray}
\begin{equation}\label{bco}
\sup_{x\in\mathbb{R}^n,\,t\in[0,\infty)}{\cal L}V_2(t,x)<\infty,
\end{equation}
and
\begin{equation}\label{b3}
\lim_{n\rightarrow\infty}\sup_{|x|>n,\,t\in[0,\infty)}{\cal
L}V_2(t,x) =-\infty.
\end{equation}
Obviously, (\textbf{H2}) implies (\textbf{H1}).

\begin{thm}\label{thm-3a}
Suppose that $(\mathbf{A1})$, $(\mathbf{A2})$, $(\mathbf{B})$ and
$(\mathbf{H^2})$ hold.
Then, the SDE (\ref{a2}) has a $\theta$-periodic solution.
\end{thm}
\noindent {\bf Proof.} Let $\{X(t),t\geq 0\}$ be the unique
solution of the SDE (\ref{a2}). For $n\in \mathbb{N}$, define
$$\beta_n=\inf\{t\in[0,\infty):|X(t)|\geq n\}.$$
For $t\ge 0$, by It\^{o}'s formula, we get
\begin{equation}\label{H43}
\mathbb{E}[V_2(t\wedge\beta_n,X(t\wedge\beta_n))]
=\mathbb{E}[V_2(0,X(0))]+\mathbb{E}\left[\int_{0}^{t\wedge\beta_n}
{\cal L}V_2(s,X(s))\mathrm{d}s\right].
\end{equation}

Define
$$
A_n:=-\sup_{|x|>n,\,t\in[0,\infty)}{\cal L}V_2(t,x).
$$
By (\ref{b3}), we get
\begin{equation}\label{H33}
\lim_{n\rightarrow\infty}A_n=\infty.
\end{equation}
We have
$$
{\cal L}V_2(t,X(s))\leq -I_{\{|X(s)|\geq
n\}}A_n+\sup_{|x|<n,t\in[0,\infty)}{\cal L}V_2(t,x).
$$
By (\ref{bq1}) and (\ref{bco}), we know that there exist positive constants $c_1$ and $c_2$ such that for large $n$,
\begin{equation}\label{H23}\mathbb{E}\left[\int_{0}^{t\wedge\beta_n}I_{\{|X(s)|\geq n\}}\mathrm{d}s\right]\leq \frac{c_{1}t+c_2}{ A_n}.\end{equation}
Denote $B_n=\{x\in\mathbb{R}^m:|x|< n\}$ and $B^c_n=\{u\in\mathbb{R}^m:|x|\ge n\}$. Letting $n\rightarrow\infty$ in (\ref{H23}), we obtain by (\ref{H33}) that
\begin{eqnarray}\label{eq1}
\lim_{n\rightarrow\infty}\varlimsup_{T\rightarrow\infty}\frac{1}{T}\int_{0}^{T}P(0,x,u,B_n^c)\mathrm{d}u=0.
\end{eqnarray}

By (\ref{bco}), there exists $\lambda>0$ such that
\begin{equation}\label{H44}{\cal L}V_2(t, x)\le \lambda,\ \ \forall t\ge 0, x\in \mathbb{R}^m.
\end{equation}
By (\ref{H43}) and  (\ref{H44}), we get
$$\mathbb{E}[V_2(t,X(t))]\leq \lambda t+V_2(0,x).$$
Together with Chebyshev's inequality, this implies that
\begin{equation}\label{H55}\mathbb{P}(0,x,t,B_n^c)\leq \frac{\lambda t+V_2(0,x)}{\inf_{|x|> n,\,t\in[0,\infty)}V_2(t,x)}.\end{equation}
By (\ref{bq1}) and (\ref{H55}), we find that there exists a sequence of positive integers $\gamma_n\uparrow\infty$ such that
\begin{eqnarray}\label{eq2}
\lim_{n\rightarrow\infty}\left\{\sup_{x\in
B_{\gamma_n},\,t\in(0,\theta)}\mathbb{P}(0,x,t,B_n^c)\right\}=0.
\end{eqnarray}

By (\ref{eq1}), (\ref{eq2}), Lemma \ref{lem-3a} and \cite[Theorem 3.2 and Remark 3.1]{a8}, we conclude that the SDE (\ref{a2}) has a $\theta$-periodic solution. Here we would like to call the reader's attention to a gap of \cite[Theorem 3.2]{a8}, which was pointed out by Hu and Xu recently. According to Hu and Xu \cite[Theorem 2.1 and Remark A.1]{a12}, \cite[Theorem 3.2 and Remark 3.1]{a8} holds under the additional assumption that $\{P_{s,t}\}$ is a Feller
semigroup. By Theorem \ref{thmFeller}, $\{P_{s,t}\}$ is a Feller semigroup and hence we can apply \cite[Theorem 3.2 and Remark 3.1]{a8} to show that the SDE (\ref{a2}) has a $\theta$-periodic solution.\hfill\fbox

\subsection{Uniqueness of periodic solutions}
\begin{lem}\label{lem-7}
Let $0\le s<t<t_1$. If a Markovian semigroup $\{P_{s,t}\}$ is strongly Feller at $(t,t_1)$ and irreducible at  $(s,t)$, then it is regular at
 $(s,t_1)$.
\end{lem}
{\bf Proof.} Suppose that $\{P_{s,t}\}$ is strongly Feller at $(t,t_1)$ and irreducible at  $(s,t)$.
Assume that for some $x_0\in\mathbb{R}^m$ and $A\in\mathcal{B}(\mathbb{R}^m)$,
$P_{s,t_1}(x_0,A)>0$. Since
\begin{eqnarray*}
P_{s,t_1}(x_0,A)&=&\int_{\mathbb{R}^m}
P(s,x_0,t,\mathrm{d}y)P(t,y,t_1,A),
\end{eqnarray*}
there exists $y_0\in \mathbb{R}^m$ such that $P(t,y_0,t_1,A)>0$. Since $P_{t,t_1}$ is strongly
Feller at  $(t,t_1)$, $P_{t,t_1}I_A\in C_b(\mathbb{R}^m)$. Hence
there exists an $r_0>0$ such that $P_{t,t_1}(y,A)>0$ for all $y\in
B(y_0,r_0)$,
where $B(y_0,r_0):=\{y\in\mathbb{R}^m: |y-y_0|<r_0\}$.
Consequently, for arbitrary $x\in\mathbb{R}^m$, we have
\begin{eqnarray*}
P_{s,t_1}(x,A)&=&\int_{\mathbb{R}^m}
P(s,x,t,\mathrm{d}y)P(t,y,t_1,A)\\
&\geq&\int_{B(y_0,r_0)} P(s,x,t,\mathrm{d}y)P(t,y,t_1,A)\\
&>&0,
\end{eqnarray*}
where we have used the fact that $P_{s,t}(x,B(y_0,r_0))>0$ and $P_{t,t_1}(y,A)>0$
for all $y\in B(y_0,r_0)$. Thus, if $P_{s,t_1}(x_0,A)>0$ for some
$x_0\in\mathbb{R}^m$ then
$P_{s,t_1}(x,A)>0$ for all $x\in\mathbb{R}^m$. Therefore, the regularity of $\{P_{s,t}\}$ at $(s,t_1)$ follows.\hfill\fbox

\begin{lem}\label{lem-8}
Let $\{P_{s,t}\}$ be a stochastically continuous $\theta$-periodic Markovian semigroup and $\{\mu_s\}$ be a family of $\theta$-periodic (probability) measures with respect to
$\{P_{s,t}\}$.
If $\{P_{s,t}\}$ is regular at $(s,s+\theta)$ for any
$s\in[0,\theta)$, then $\{\mu_s\}$ is the unique $\theta$-periodic
measures with respect to $\{P_{s,t}\}$.
\end{lem}
{\bf Proof.} Let $\{\mu_s\}$ be a family of $\theta$-periodic   measures with respect to
$\{P_{s,t}\}$. Suppose that $\{P_{s,t}\}$ is regular at $(s,s+\theta)$ for any
$s\in[0,\theta)$.

\emph{Step 1.} We first show that $\mu_s$ is ergodic for any $s\ge
0$. That is,  if $A\in\mathcal{B}(\mathbb{R}^m)$ and
\begin{equation}\label{va}
  P(s,x,s+\theta,A)
=I_A(x),\ \ \ \ \mu_s- a.s.,
\end{equation}
then $\mu_s(A)=0$ or $\mu_s(A)=1$.

Let $A\in\mathcal{B}(\mathbb{R}^m)$ satisfying
$\mu_s(A)>0$. We will show that $\mu_s(A)=1$.
By (\ref{va}), we get
$$\mu_s (\{x\in A:P(s,x,s+\theta,A)=1\})=\mu_s(A).$$
Then, there exists $x_0\in A$ such that
$P(s,x_0,s+\theta,A)=1$.
Since all probabilities $P(s,x,s+\theta,\cdot)$, $x\in \mathbb{R}^m$, are mutually equivalent,
we get $P(s,x,s+\theta,A)=P(s,x_0,s+\theta,A)=1$, $\forall
x\in\mathbb{R}^m$.
It follows that
$$\mu_s(A)=\int_{\mathbb{R}^m}P(s,y,s+\theta,A)\mu_s(dy)=1.$$

\emph{Step 2.} Next we show that for any
$s\ge 0$ and $\varphi\in L^2(\mathbb{R}^m;\mu_s)$,
\begin{equation}\label{a1}
\lim_{n\rightarrow\infty}
\frac{1}{n}\sum_{i=1}^{n}P_{s,s+i\theta}\varphi=\int_{\mathbb{R}^m}\varphi
d\mu_s\ \
{\rm in}\ \ L^2(\mathbb{R}^m;\mu_s).
\end{equation}
Note that
\begin{eqnarray*}
\left\|\frac{1}{n}\sum_{i=1}^{n}P_{s,s+i\theta}\varphi\right\|_{L^2(\mathbb{R}^m;\mu_s)}&\le&\frac{1}{n}\sum_{i=1}^{n}\|P_{s,s+i\theta}\varphi\|_{L^2(\mathbb{R}^m;\mu_s)}\\
&\le&\frac{1}{n}\sum_{i=1}^{n}\left(\int_{\mathbb{R}^m}P_{s,s+i\theta}(\varphi^2)
d\mu_s\right)^{1/2}\\
&=&\left(\int_{\mathbb{R}^m}\varphi^2 d\mu_s\right)^{1/2}.
\end{eqnarray*}
Hence, in order to prove (\ref{a1}), it is sufficient to show that
\begin{equation}\label{wa2}
\lim_{n\rightarrow\infty}
\frac{1}{n}\sum_{i=1}^{n}P_{s,s+i\theta}\varphi=\int_{\mathbb{R}^m}\varphi
d\mu_s\ \
{\rm weakly\ in}\ \ L^2(\mathbb{R}^m;\mu_s).
\end{equation}

Define $\Xi:=\mathbb{R}^{\mathbb{Z}}$. Let ${\cal G}$ be the
$\sigma$-algebra generated by the set of all cylindrical sets on $\Xi$. By
the Kolmogorov extension theorem, there exists a unique probability measure
$\mathbb{P}_{\mu_s}$ on $(\Xi,{\cal G})$ such that
\begin{eqnarray*}
&&\mathbb{P}_{\mu_s}(\{\omega\in\Xi:\omega_{n_1}\in A_1,\omega_{n_2}\in
A_2,\dots,\omega_{n_k}\in A_k,n_1<n_2<\cdots<n_k\})\\
&&=\int_{A_1}\mu_s(dx_1)\int_{A_2}
P(s, x_1,s+\theta,dx_2)\dots\int_{A_k}P(s,x_{k-1},s+(k-1)\theta,dx_k),\\
&&\ \ \ \ \ \ \ \ \ \ \ \ \forall A_i\in {\cal B}(\mathbb{R}^m),\ 1\le i\le
k,\ k\in\mathbb{N}.
\end{eqnarray*}
Define $\Theta:\Xi\rightarrow\Xi$ by $(\Theta\omega)_l=\omega_{l+1}$,
$l\in\mathbb{Z}$. Then, $\Theta$ is a measure preserving transformation.
Let $\varphi\in L^2(\mathbb{R}^m;\mu_s)$ and define
\begin{equation}\label{lkl}
\xi=\varphi(\omega_0).
 \end{equation}
We have $\xi\in L^2(\Omega,{\cal F},\mathbb{P}_{\mu_s})$. Then, we obtain by
Birkhoff's ergodic theorem (cf. \cite[Theorem 2.3]{P125}) that there exits
$\xi^*\in L^2(\Xi,{\cal G},\mathbb{P}_{\mu_s})$ such that
\begin{equation}\label{wan1}
\lim_{n\rightarrow\infty}
\frac{1}{n}\sum_{i=1}^{n}\xi(\Theta^i)=\xi^*,\ \ P_{\mu_s}\mbox{-}a.s.\ {\rm
and\ in}\ \ L^2(\Xi,{\cal G},\mathbb{P}_{\mu_s}).
\end{equation}
Define
\begin{equation}\label{wan2}
U_i\eta(\omega)=\eta(\Theta^i\omega),\ \ \eta\in L^2(\Xi,{\cal
G},\mathbb{P}_{\mu_s}),\ \ \omega\in\Xi,\ i\in\mathbb{N}.
\end{equation}
Then, we obtain by (\ref{wan1}) that
\begin{equation}\label{wan3}
U_i\xi^*=\xi^*,\ \ i\in \mathbb{N}.
\end{equation}

Define ${\cal G}_j=\sigma\{\omega_u:u\le j\}$, $j\in\mathbb{Z}$, and ${\cal
G}_{[r,j]}=\sigma\{\omega_u:r\le u\le j\}$ for $r\le j\in\mathbb{Z}$.
Following the argument of the proof of \cite[Lemma 3.2.2]{3a}, we can show
that for any $\eta$ which is ${\cal G}_{[-j,j]}$-measurable, $j\in
\mathbb{N}$,
$$
\mathbb{E}_{{\mu_s}}[|\mathbb{E}_{{\mu_s}}(U_j\eta|{\cal
G}_{[0,0]})-\xi^*|^2]\le 10E_{{\mu_s}}[|\eta-\xi^*|^2].
$$
Following the argument of the proof of \cite[Propositon 2.2.1]{3a}, we can
show that for arbitrary $F\in{\cal G}$ and $\varepsilon>0$ there exists a
cylindrical set $C$ such that
$$
\mathbb{P}_{\mu_s}(F\backslash C)+\mathbb{P}_{\mu_s}(C\backslash
F)<\varepsilon.
$$
Then, there exists a sequence $\{\eta_j\}$ of ${\cal G}_{[-j,j]}$-measurable
elements of $L^2(\Xi,{\cal G},P_{\mu_s})$ such that
$$
\lim_{j\rightarrow\infty}\mathbb{E}_{{\mu_s}}[U_j\eta_j|{\cal
G}_{[0,0]}]=\xi^*\ \ {\rm in}\  L^2(\Xi,{\cal G},P_{\mu_s}).
$$
Moreover, there exists $\{\varphi_j\}\subset L^2(\Xi,{\cal
G},\mathbb{P}_{\mu_s})$ such that
$$
\mathbb{E}_{{\mu_s}}[U_j\eta_j|{\cal G}_{[0,0]}]=\varphi_j(\omega_0),\ \
P_{\mu_s}\mbox{-}a.s..
$$
Without loss of generality we can assume that
$$
\lim_{j\rightarrow\infty}\varphi_j(\omega_0)=\xi^*,\ \
\mathbb{P}_{\mu_s}\mbox{-}a.s.\ {\rm and\ in}\ L^2(\Xi,{\cal
G},P_{\mu_s}).
$$

Define
$$
\varphi(x)=\left\{\begin{array}{ll}
                  \lim_{j\rightarrow\infty}\varphi_j(x),\ \ &{\rm if\ the\
                  limit\ exists},\\
                  0, &{\rm otherwise}.
                \end{array}\right.
$$
Then, we have $\xi^*(\omega)=\varphi(\omega_0)$,
$\mathbb{P}_{\mu_s}\mbox{-}a.s.$. By (\ref{wan2}) and
(\ref{wan3}), we get
\begin{equation}\label{wan4}
\varphi(\omega_0)=\xi^*(\omega)=\xi^*(\Theta^i\omega)=\varphi((\Theta^i\omega)_0)=\varphi(\omega_i),\
\ i\in \mathbb{N}.
\end{equation}
We claim that $\varphi$ is a constant. In fact, define
$\Lambda=\varphi^{-1}(\alpha,\infty)$ for $\alpha\in \mathbb{R}$.
Then, we obtain by (\ref{wan4}) that
$$
P(s,x,s+\theta,\Lambda)=\mathbb{E}_{s,x}[\chi_{\Lambda}(\omega_1)]=\mathbb{E}_{s,x}[\chi_{\Lambda}(\omega_0)]=\chi_\Lambda(x),\
\ \mu_s-a.s..
$$
Since $\mu_s$ is ergodic by Step 1, we have that $\mu(\Lambda)=0$ or $1$. Since
$\alpha\in \mathbb{R}$ is arbitrary, we conclude that $\varphi$ is
a constant. Thus, $\xi^*$ is a constant. Therefore, we obtain by
(\ref{lkl}) and (\ref{wan1}) that
\begin{equation}\label{lkl1}\lim_{n\rightarrow\infty}
\frac{1}{n}\sum_{i=1}^{n}\varphi(\omega(i))=\int_{\mathbb{R}^m}\varphi d\mu_s\
\ {\rm in}\ \ L^2(\Xi,{\cal G},P_{\mu_s}).
\end{equation}
Let $\zeta\in L^2(\mathbb{R}^m;\mu_s)$. Then, we obtain by (\ref{lkl1}) that
$$\lim_{n\rightarrow\infty}
\frac{1}{n}\sum_{i=1}^{n}\langle
P_{s,s+i\theta}\phi,\zeta\rangle_{L^2(\mathbb{R}^m;\mu_s)}=\int_{\mathbb{R}^m}\varphi
d\mu_s\int_{\mathbb{R}^m}\zeta d\mu_s.
$$
Since $\zeta\in L^2(\mathbb{R}^m;\mu_s)$ is arbitrary, the proof of
(\ref{wa2}) is complete.

\emph{Step 3.} Finally, we show that $\{\mu_s\}$ is the unique $\theta$-periodic measures. Suppose
that $\{\mu'_s\}$ is a different family of $\theta$-periodic measures with respect to
$\{P_{s,t}\}$. Then, there exist $s\ge 0$ and  $\Gamma\in {\cal B}(\mathbb{R}^m)$
such that
$$
\mu_s(\Gamma)\not=\mu'_s(\Gamma).
$$

By Step 2,  there exists a sequence $\{T_n\uparrow\infty\}$ such that
$$\lim_{n\rightarrow\infty}\frac{1}{T_n}\sum_{i=1}^{T_n}p(s,x,s+i\theta,\Gamma)=\mu_s(\Gamma),\
\ \mu_s\mbox{-}a.s.,
$$
and
$$
\lim_{n\rightarrow\infty}\frac{1}{T_n}\sum_{i=1}^{T_n}p(s,x,s+i\theta,\Gamma)=\mu'_s(\Gamma),\
\ \mu'_s\mbox{-}a.s..$$
Define
$$\left\{x\in\mathbb{R}^m:\lim_{n\rightarrow\infty}\frac{1}{T_n}
\sum_{i=1}^{T_n}p(s,x,s+\theta,\Gamma)=\mu_s(\Gamma)\right\}=A,$$
and
$$\left\{x\in\mathbb{R}^m:\lim_{n\rightarrow\infty}\frac{1}{T_n}
\sum_{i=1}^{T_n}p(s,x,s+\theta,\Gamma)=\mu'_s(\Gamma)\right\}=B.$$
It is clear that $A\cap B=\emptyset$ and
$\mu_s(A)=\mu'_s(B)=1$. Thus, $\mu_s$ and $\mu'_s$ are singular.
However, since $\{P_{s,t}\}$ is regular at $(s,s+\theta)$ and $\{\mu_s\}$ and
$\{\mu'_s\}$ are $\theta$-periodic measures (cf. (\ref{perio})), we must have that $\mu_s$ and
$\mu'_s$ are equivalent to $P(s,x,s+\theta,\cdot)$, $x\in
\mathbb{R}^m$. We have arrived at a contradiction. \hfill\fbox

We put the following assumption, which implies the conditions  $\mathbf{(\mathbf{H^2})}$ and $(\mathbf{H^3})$.

\noindent $(\mathbf{H})$\ \  (i) There exists $W(t,x):[0,\infty)\times\mathbb{R}^m\rightarrow\mathbb{R}^m$ which is locally bounded and for each $n\in\mathbb{N}$, there exists
$R_n\in L_{loc}^1([0,\infty);\mathbb{R}_+)$ such that for any $t\in
[0,\infty)$ and $x,y\in\mathbb{R}^m$ with $|x|\vee|y|\leq n$,
$$
   |W(t,x)-W(t,y)|^2\leq
  R_n(t)|x-y|^2.
  $$

  \noindent (ii) There exist $q\in B_{b,loc}(\mathbb{R}_+)$ and
$U(t,x):\mathbb{R}_+\times \mathbb{R}_+\rightarrow(0,\infty)$ which is increasing with respect to $x$.

  \noindent (iii) There exists $V\in
C^{1,2}([0,\infty)\times\mathbb{R}^m,\mathbb{R}_+)$  such that
$$ \lim_{|x|\rightarrow\infty}
\left[\inf_{t\in[0,\infty)}V(t,x)\right]=\infty,
$$
$$
\sup_{x\in\mathbb{R}^n,\,t\in[0,\infty)}{\cal L}V(t,x)<\infty,
$$
$$
\lim_{n\rightarrow\infty}\sup_{|x|>n,\,t\in[0,\infty)}{\cal
L}V(t,x) =-\infty,
$$
and for $t\ge 0$ and $x\in\mathbb{R}^m$,
$$
U(t,|x|)\le V(t,x)\le\langle W,V_x(t,x)\rangle+q(t).
$$

\begin{thm}\label{lv}
Suppose that $\mathbf{(A3)}$, $\mathbf{(A4)}$ and $(\mathbf{H})$ hold.
Then,

(i) The SDE (\ref{a2}) has a unique $\theta$-periodic solution $\{X(t),t\geq 0\}$;

 (ii) The Markovian transition semigroup $\{P_{s,t}\}$ of $\{X(t),t\geq 0\}$ is strongly Feller and irreducible;

 (iii) Let $\mu_s(A)=\mathbb{P}(X(s)\in A)$ for $A\in
\mathcal{B}(\mathbb{R}^m)$ and $s\ge 0$. Then,
for any
$s\ge 0$ and $\varphi\in L^2(\mathbb{R}^m;\mu_s)$, we have
$$
\lim_{n\rightarrow\infty}
\frac{1}{n}\sum_{i=1}^{n}P_{s,s+i\theta}\varphi=\int_{\mathbb{R}^m}\varphi
d\mu_s\ \
{\rm in}\ \ L^2(\mathbb{R}^m;\mu_s).
$$
\end{thm}
{\bf Proof.} By $(\mathbf{A3})$, $(\mathbf{H})$ and
Theorem \ref{lem-3}, we know that the SDE (\ref{a2}) has a unique solution. Further, by $(\mathbf{A3})$, $(\mathbf{A4})$,
$(\mathbf{H})$ and Theorem \ref{thm-s}, we know that  $\{P_{s,t}\}$ is strongly Feller. Hence the SDE (\ref{a2}) has a $\theta$-periodic solution by
Theorem \ref{thm-3a}. The uniqueness of the $\theta$-periodic solution is a direct consequence of Theorem \ref{thm-s}, Theorem \ref{thm-6}, Lemma \ref{lem-7} and
Lemma \ref{lem-8}. Finally, the last assertion of the theorem follows from the proof of Lemma \ref{lem-8} (see (\ref{a1})). \hfill\fbox

\section{Examples}\label{sec4}\setcounter{equation}{0}
In this section, we use three examples to illustrate the main
results of the paper.

\begin{exa} Suppose that the coefficient functions  $b,\sigma,H,G$ of the SDE (\ref{a2}) are $\theta$-periodic, $b(\cdot,0)\in
{L}^2([0,\theta);\mathbb{R}^m)$, and for each $n\in\mathbb{N}$, there exists
$L_n\in L^{\infty}([0,\theta);\mathbb{R}_+)$ such that for any $t\in
[0,\theta)$ and $x,y\in\mathbb{R}^m$ with $|x|\vee|y|\leq n$,
\begin{eqnarray*}
 &&   | b(t,x)-b(t,y)|^2\leq
  L_n(t)|x-y|^2,\ \ \ \
|\sigma(t,x)-\sigma(t,y)|^2\leq
  L_n(t)|x-y|^2,  \\
 &&   \int_{\{|u|<1\}}|H(t,x,u)-H(t,y,u)|^2\nu(\mathrm{d}u)\leq
L_n(t)|x-y|^2.
  \end{eqnarray*}
For any $t\in[0,\theta)$ and $x\in\mathbb{R}^m$,
$Q(t,x)=\sigma(t,x)\sigma^T(t,x)$ is invertible and
$$
\sup_{|x|\leq
n,\,t\in[0,\theta)}| Q^{-1}(t,x)|<\infty,\ \ \forall n\in \mathbb{N}.
$$

In addition, we assume that there exist $r,c_1,c_2>0$ such that
$$
\langle b(t,x),x\rangle+|\sigma(t,x)|^2\leq-c_1|x|^{r}+c_2,\ \ \forall t\in[0,\theta),x\in\mathbb{R}^m,
$$
and for any $\varepsilon>0$ there exists $c_{\varepsilon}>0$ such that
\begin{eqnarray*}
&&\int_{\{|u|<1\}}|H(t,x,u)|^2\nu(du)+|x|\int_{\{|u|\ge 1\}}|G(t,x,u)|\nu(du)\\
&&\ \ \ \ \ \ +\int_{\{|u|\ge 1\}}|G(t,x,u)|^2\nu(du)\le\varepsilon|x|^{r}+c_{\varepsilon},\ \ \forall t\in[0,\theta),x\in\mathbb{R}^m.
\end{eqnarray*}

Define $$
V(t,x)=|x|^2+1,
$$
and
$$W(t,x)=\frac{x}{2}.$$
Then,
\begin{eqnarray*}
   V(t,x)=\langle W, V_x(t,x)\rangle+1.
\end{eqnarray*}
We have
\begin{eqnarray*}
&&{\int_{\{|u|<1\}}}[V(t,x+H(t,x,u))-V(t,x)-\langle V_x(t,x),H(t,x,u)\rangle]\nu(\mathrm{d}u)\\
&=&
{\int_{\{|u|<1\}}}\frac{1}{2}H^T(t,x,u)
    V_{xx}(x+\theta_1H(t,x,u))H(t,x,u)\nu(\mathrm{d}u)\\
&\leq&\int_{\{|u|<1\}}|H(t,x,u)|^2\nu(du),
\end{eqnarray*}
and
\begin{eqnarray*}
&&{\int_{\{|u|\ge 1\}}}[V(t,x+G(t,x,u))-V(t,x)]\nu(\mathrm{d}u)\\
&=&{\int_{\{|u|\ge 1\}}}\left[2\langle x,G(t,x,u)\rangle +|G(t,x,u)|^2\right]\nu(\mathrm{d}u)\\
&\le&2|x|\int_{\{|u|\ge 1\}}|G(t,x,u)|\nu(du)+\int_{\{|u|\ge 1\}}|G(t,x,u)|^2\nu(du).
\end{eqnarray*}
Hence, for any $\varepsilon>0$, we have
\begin{eqnarray*}
    {\cal L}V(t,x)&=&\langle V_x(t,x), b(t,x)\rangle+\frac{1}{2}\mathrm{ trace}(\sigma^T(t,x)V_{xx}(t,x)\sigma(t,x))\\
    &&+{\int_{\{|u|<1\}}}[V(t,x+H(t,x,u))-V(t,x)
    -\langle V_x(t,x), H(t,x,u)\rangle]\nu(\mathrm{d}u)\\
    &&+{\int_{\{|u|\ge 1\}}}[V(t,x+G(t,x,u))-V(t,x)]\nu(\mathrm{d}u)\\
    &\leq&-c_1|x|^{r}+c_2+2(\varepsilon|x|^{r}+c_{\varepsilon}).
  \end{eqnarray*}
Then, $$\sup_{x\in\mathbb{R}^n,\,t\in[0,\infty)}{\cal L}V(t,x)<\infty,
$$
 and
 $$\lim_{n\rightarrow\infty}\sup_{|x|>n,\,t\in[0,\infty)}{\cal
L}V(t,x) =-\infty.
$$
 Thus all conditions of Theorem \ref{lv} are satisfied and therefore all assertions of Theorem \ref{lv} hold.
\end{exa}

\begin{exa} (Stochastic Lorenz equation)

The Lorenz equation is a remarkable mathematical model
for atmospheric convection, which was introduced by E.N. Lorenz in \cite{lo}. In recent years, many papers have been devoted to the Lorenz equation with noises (cf.  \cite{Lo1} and the references therein).

 In this example, we consider the following Lorenz equation with multiplicative L\'evy noise:
\begin{eqnarray}\label{Lorenz}
dX_1(t) &=& (-\alpha(t)X_1(t-)+\alpha(t)X_2(t-))dt+\sum_{j=1}^3\sigma_{1j}(t,X(t-))dB_{j}(t)\nonumber\\
&& +\int_{\{|u|<1\}}H_1(t,X(t-),u)\widetilde{N}(dt,du)+\int_{\{|u|\ge 1\}}G_1(t,X(t-),u)N(dt,du),\nonumber\\
  dX_2(t) &=& (\mu(t)X_1(t-)-X_2(t-)-X_1(t-)X_3(t-))dt+\sum_{j=1}^3\sigma_{2j}(t,X(t-))dB_{j}(t)\nonumber\\
\nonumber\\
 &&+\int_{\{|u|<1\}}H_2(t,X(t-),u)\widetilde{N}(dt,du)
 +\int_{\{|u|\ge 1\}}G_2(t,X(t-),u)N(dt,du),\nonumber\\
 dX_3(t) &=& (-\beta(t)X_3(t-)+X_1(t-)X_2(t-))dt+\sum_{j=1}^3\sigma_{3j}(t,X(t-))dB_{j}(t)\nonumber\\
 &&+\int_{\{|u|<1\}}H_3(t,X(t-),u)\widetilde{N}(dt,du)+\int_{\{|u|\ge 1\}}G_3(t,X(t-),u)N(dt,du),
\end{eqnarray}
where
$\alpha(t),\beta(t),\mu(t):[0,\infty)\rightarrow\mathbb{R}_+$,
$\sigma(t,x):[0,\infty)\times\mathbb{R}^3\rightarrow
\mathbb{R}^{3\times 3}$,
$H(t,x,u):[0,\infty)\times\mathbb{R}^3\times\mathbb{R}^l
\rightarrow\mathbb{R}^3$ and
$G(t,x,u):[0,\infty)\times\mathbb{R}^3\times\mathbb{R}^l
\rightarrow\mathbb{R}^3$ are all Borel measurable and
$\theta$-periodic.

We assume that $\alpha(t),\beta(t),\mu(t)$ are continuously differentiable with
$$
\min\{\alpha(t):t\in [0,\theta)\}>0,\ \ \ \ \min\{\beta(t):t\in [0,\theta)\}>0.
$$
For each $n\in\mathbb{N}$, there exists
$L_n\in L^{\infty}([0,\theta);\mathbb{R}_+)$ such that for any $t\in
[0,\theta)$ and $x,y\in\mathbb{R}^m$ with $|x|\vee|y|\leq n$,
\begin{eqnarray*}
 &&
|\sigma(t,x)-\sigma(t,y)|^2\leq
  L_n(t)|x-y|^2,  \\
 &&   \int_{\{|u|<1\}}|H(t,x,u)-H(t,y,u)|^2\nu(\mathrm{d}u)\leq
L_n(t)|x-y|^2.
  \end{eqnarray*}
For any $t\in[0,\theta)$ and $x\in\mathbb{R}^m$,
$Q(t,x)=\sigma(t,x)\sigma^T(t,x)$ is invertible and
$$
\sup_{|x|\leq
n,\,t\in[0,\theta)}| Q^{-1}(t,x)|<\infty,\ \ \forall n\in \mathbb{N}.
$$
Moreover, for any $\varepsilon>0$ there exists $c_{\varepsilon}>0$ such that
\begin{eqnarray*}
&&|\sigma(t,x)|^2+\int_{\{|u|<1\}}|H(t,x,u)|^2\nu(du)\\
&&\ \ \ \ \ \ +\int_{\{|u|\ge 1\}}|G(t,x,u)|^2\nu(du)
\leq \varepsilon|x|^{2}+c_{\varepsilon},\ \ \forall t\in[0,\theta),x\in\mathbb{R}^m.
\end{eqnarray*}

Define
$$V(t,x)=x_1^2+x_2^2+(x_3-\alpha(t)-\mu(t))^2+(\alpha(t)+\mu(t))^2+1.
$$
We have that
$$\frac{\partial V}{\partial x_1}=2x_1,\ \
\frac{\partial V}{\partial x_2}=2x_2,\ \
\frac{\partial V}{\partial x_3}=2(x_3-\alpha(t)-\mu(t)).$$
Let
$$
U(t,x)=\frac{|x|^2}{2}+1,
$$
and
$$W(t,x)=\frac{1}{2}(x_1,x_2,x_3-\alpha(t)-\mu(t)).$$
Then,
\begin{eqnarray*}
   U(t,|x|)\le V(t,x)=\langle W,V_x(t,x)\rangle+(\alpha(t)+\mu(t))^2+1.
\end{eqnarray*}

Hence, for any $\varepsilon>0$, we have
\begin{eqnarray*}
  {\cal L}V(t,x) &=& 2(2\alpha(t)+2\mu(t)-x_3)(\alpha'(t)+\mu'(t))\\
  &&-2\alpha(t)x_1^2-2x_2^2-2\beta(t) (x_3^2-(\alpha(t)+\mu(t))x_3)+|\sigma(t,x)|^2\\
  &&+\sum_{i=1}^3\int_{\{|u|<1\}}|H_i(t,x,u)|^2\nu(du)
   +\sum_{i=1}^2\int_{\{|u|\ge 1\}}[(x_i+G_i(t,x,u))^2-x_i^2]\nu(du)\\
   &&+\int_{\{|u|\ge 1\}}[(x_3+G_3(t,x,u)-\alpha(t)-\mu(t))^2-(x_3-\alpha(t)-\mu(t))^2]\nu(du)\\
 &\leq&2(2\alpha(t)+2\mu(t)+|x_3|)(|\alpha'(t)|+|\mu'(t)|)\\
 &&-2\left[\alpha(t)x_1^2+x_2^2+\beta(t)x_3^2-\beta(t)(\alpha(t)+\mu(t))x_3\right]+\varepsilon|x|^{2}+c_{\varepsilon}\\
 &&+2[|x_1|+|x_2|+|x_3|+\alpha(t)+\mu(t)]\left[\nu(\{|u|\ge1\})\right]^{1/2}(\varepsilon|x|^{2}+c_{\varepsilon})^{1/2}.
\end{eqnarray*}
Then, $$\sup_{x\in\mathbb{R}^n,\,t\in[0,\infty)}{\cal L}V(t,x)<\infty,
$$
 and
 $$\lim_{n\rightarrow\infty}\sup_{|x|>n,\,t\in[0,\infty)}{\cal
L}V(t,x) =-\infty.
$$
 Thus all conditions of Theorem \ref{lv} are satisfied.
Therefore, the stochastic Lorenz equation (\ref{Lorenz}) has a unique $\theta$-periodic solution $\{X(t),t\geq 0\}$ and assertions (ii) and (iii) of Theorem \ref{lv} hold.
\end{exa}

\begin{exa} (Equation of the lemniscate of Bernoulli with L\'evy noise)

In this example, we consider the stochastic equation of the lemniscate of Bernoulli, which generalizes \cite[Example 3]{CD} to the non-autonomous case with L\'evy noise.

For $x=(x_1,x_2)\in \mathbb{R}^2$, define
$$I(x)=(x_1^2+x_2^2)^2-4(x_1^2-x_2^2).
$$ Let
$$
{\cal V}(I)=\frac{I^2}{2(1+I^2)^{3/4}},\ \ \ \ {\cal H}(I)=\frac{I}{(1+I^2)^{3/8}}.
$$
Consider the vector field
$$
b(x)=-\left[{\cal V}_x(I)+\left(\frac{\partial {\cal H}(I)}{\partial x_2},-\frac{\partial {\cal H}(I)}{\partial x_1}\right)^T\right].
$$
We have
\begin{eqnarray}
&&\frac{d{\cal V}(I)}{dI}=\frac{I(I^2+4)}{4(1+I^2)^{7/4}},\ \ \ \ \frac{d{\cal H}(I)}{dI}=\frac{I^2+4}{4(1+I^2)^{11/4}},\label{lemirr}\\
&&\frac{\partial I}{\partial x_1}=4x_1(x_1^2+x_2^2)-8x_1,\ \ \ \ \frac{\partial I}{\partial x_2}=4x_2(x_1^2+x_2^2)+8x_2.\label{lemi1}
\end{eqnarray}
Define
\begin{eqnarray}\label{lemi2}
f(I)=\frac{d{\cal V}(I)}{dI},\ \ \ \ g(I)=\frac{d{\cal H}(I)}{dI}.
\end{eqnarray}
Then,
\begin{eqnarray}\label{lemi3}{\cal V}_x(I)=\frac{d{\cal V}(I)}{dI}\left(\frac{\partial I}{\partial x_1},\frac{\partial I}{\partial x_2}\right)^T,
\end{eqnarray}
and
\begin{eqnarray}\label{lemi4}
&&b_1(x)=-f(I)(4x_1(x_1^2+x_2^2)-8x_1)-g(I)(4x_2(x_1^2+x_2^2)+8x_2),\nonumber\\
&&b_2(x)=-f(I)(4x_2(x_1^2+x_2^2)+8x_2)-g(I)(-4x_1(x_1^2+x_2^2)+8x_1).
\end{eqnarray}

We consider the following SDE:
\begin{eqnarray}\label{ssll}
dX_1(t) &=& b_1(X(t-))dt+\sigma_{11}(t,X(t-))dB_1(t)+\sigma_{12}(t,X(t-))dB_2(t)\nonumber\\
 &&+\int_{\{|u|<1\}}H_1(t,X(t-),u)\widetilde{N}(dt,du)+\int_{\{|u|\ge 1\}}G_1(t,x(t-),u)N(dt,du),\nonumber\\
  dX_2(t) &=& b_2(X(t-))dt+\sigma_{21}(t,X(t-))dB_1(t)+\sigma_{22}(t,X(t-))dB_2(t)\nonumber\\
 &&+\int_{\{|u|<1\}}H_2(t,X(t-),u)\widetilde{N}(dt,du)
 +\int_{\{|u|\ge 1\}}G_2(t,x(t-),u)N(dt,du),
\end{eqnarray}
where $\sigma(t,x):[0,\infty)\times\mathbb{R}^2\rightarrow
\mathbb{R}^{2\times 2}$,
$H(t,x,u):[0,\infty)\times\mathbb{R}^2\times\mathbb{R}^l
\rightarrow\mathbb{R}^2$ and
$G(t,x,u):[0,\infty)\times\mathbb{R}^2\times\mathbb{R}^l
\rightarrow\mathbb{R}^2$ are all Borel measurable and  $\theta$-periodic.

We assume that for each $n\in\mathbb{N}$, there exists
$L_n\in L^{\infty}([0,\theta);\mathbb{R}_+)$ such that for any $t\in
[0,\theta)$ and $x,y\in\mathbb{R}^m$ with $|x|\vee|y|\leq n$,
\begin{eqnarray*}
 &&
|\sigma(t,x)-\sigma(t,y)|^2\leq
  L_n(t)|x-y|^2,  \\
 &&   \int_{\{|u|<1\}}|H(t,x,u)-H(t,y,u)|^2\nu(\mathrm{d}u)\leq
L_n(t)|x-y|^2.
  \end{eqnarray*}
For any $t\in[0,\theta)$ and $x\in\mathbb{R}^m$,
$Q(t,x)=\sigma(t,x)\sigma^T(t,x)$ is invertible and
$$
\sup_{|x|\leq
n,\,t\in[0,\theta)}| Q^{-1}(t,x)|<\infty,\ \ \forall n\in \mathbb{N}.
$$
Moreover, for any $\varepsilon>0$ there exists $c_{\varepsilon}>0$ such that
\begin{eqnarray*}
&&|\sigma(t,x)|^2+\int_{\{|u|<1\}}|H(t,x,u)|^2\nu(du)+|x|\int_{\{|u|\ge 1\}}|G(t,x,u)|\nu(du)\\
&&\ \ \ \ \ \ +\int_{\{|u|\ge 1\}}|G(t,x,u)|^2\nu(du)
\leq \varepsilon|x|^{2}+c_{\varepsilon},\ \ \forall t\in[0,\theta),x\in\mathbb{R}^m.
\end{eqnarray*}

Define
$$V(t,x)={\cal V}(I(x))+2^{6},\ \ \ \ t\ge0,\,x\in \mathbb{R}^2.
$$ We have that
\begin{equation}\label{estimate1}
I(x)\ge \frac{x_1^4+2x_1^2x_2^2+ x_2^4}{2}+\frac{x_1^4}{2}-4x_1^2\ge \frac{|x|^4}{2}\ \ \ \ {\rm for}\ |x|\ge 16\ {\rm with}\ |x_1|\ge |x_2|.
\end{equation}
Let
$$U(t,x)=\frac{|x|^2}{2^{9/4}}+1.$$
By (\ref{estimate1}), we get
$$
U(t,|x|)\le V(t,x),\ \ \forall t\ge0,\,x\in \mathbb{R}^2.
$$
Let
$$W(t,x)=x.$$
Then,
\begin{eqnarray*}
V(t,x)&\le&\frac{2I^2+\frac{I^4}{2}}{(1+I^2)^{7/4}}+2^{6}\\
&\le&\frac{(I+\frac{I^3}{4})(4(x_1^2+x_2^2)^2-8(x_1^2-x_2^2))}{(1+I^2)^{7/4}}+2^{6}\\
&=&\langle W,V_x(t,x)\rangle+2^{6},\ \ \forall t\ge0,\,x\in \mathbb{R}^2.
\end{eqnarray*}

By (\ref{lemi1}), (\ref{lemi2}) and (\ref{lemi4}), we get
\begin{eqnarray}\label{llemi1}
\varlimsup_{|x|\rightarrow\infty}\frac{\langle V_x(t,x), b(x)\rangle}{|x|^2}<0.
\end{eqnarray}
By direct calculation, we find that there exists a constant $\varrho_1>0$ such that
\begin{eqnarray}\label{llemi0}
\left|\frac{\partial ^2V(x)}{\partial x_1^2}\right|,\, \left|\frac{\partial ^2V(x)}{\partial x_2^2}\right|,\, \left|\frac{\partial ^2V(x)}{\partial x_1\partial x_2}\right|\le \varrho_1,\ \ \ \ \forall x\in \mathbb{R}^2.
\end{eqnarray}
Then,
\begin{eqnarray}\label{llemi2}
&&{\int_{\{|u|<1\}}}[V(t,x+H(t,x,u))-V(t,x)-\langle V_x(t,x),H(t,x,u)\rangle]\nu(\mathrm{d}u)\nonumber\\
&=&
{\int_{\{|u|<1\}}}\frac{1}{2}H^T(t,x,u)
    V_{xx}(x+\theta_1H(t,x,u))H(t,x,u)\nu(\mathrm{d}u)\nonumber\\
&\leq&\varrho_1\int_{\{|u|<1\}}|H(t,x,u)|^2\nu(du),
\end{eqnarray}
where $\theta_1\in (0,1)$.
By (\ref{lemirr}), (\ref{lemi1}) and (\ref{lemi3}), we find that there exists  a constant $\varrho_2>0$, which is independent of $t,x$, such that
\begin{eqnarray}\label{llemi3}
&&{\int_{\{|u|\ge 1\}}}[V(t,x+G(t,x,u))-V(t,x)]\nu(\mathrm{d}u)\nonumber\\
&=&{\int_{\{|u|\ge 1\}}}V_{x}(x+\theta_2G(t,x,u))G(t,x,u)\nu(\mathrm{d}u)\nonumber\\
&\le&\varrho_2{\int_{\{|u|\ge 1\}}}(1+|x|+|G(t,x,u)|)|G(t,x,u)|\nu(\mathrm{d}u),
\end{eqnarray}
where $\theta_2\in (0,1)$.

Note that
\begin{eqnarray*}
    {\cal L}V(t,x)&=&\langle V_x(t,x), b(t,x)\rangle+\frac{1}{2}\mathrm{ trace}(\sigma^T(t,x)V_{xx}(t,x)\sigma(t,x))\\
    &&+{\int_{\{|u|<1\}}}[V(t,x+H(t,x,u))-V(t,x)
    -\langle V_x(t,x), H(t,x,u)\rangle]\nu(\mathrm{d}u)\\
    &&+{\int_{\{|u|\ge 1\}}}[V(t,x+G(t,x,u))-V(t,x)]\nu(\mathrm{d}u).
  \end{eqnarray*}
Then, we obtain by (\ref{llemi1})--(\ref{llemi3}) that
$$\sup_{x\in\mathbb{R}^n,\,t\in[0,\infty)}{\cal L}V(t,x)<\infty,
$$
 and
 $$\lim_{n\rightarrow\infty}\sup_{|x|>n,\,t\in[0,\infty)}{\cal
L}V(t,x) =-\infty.
$$
 Thus all conditions of Theorem \ref{lv} are satisfied.
Therefore, the stochastic equation of the lemniscate of Bernoulli (\ref{ssll}) has a unique $\theta$-periodic solution  $\{X(t),t\geq 0\}$ and assertions (ii) and (iii) of Theorem \ref{lv} hold.
\end{exa}

\bigskip

{ \noindent {\bf\large Acknowledgments}\ \
This work was supported by the Graduate Joint Training Program of the Guangdong Educational Department, China, and
the  Natural Sciences and Engineering Research Council of Canada.


\begin{thebibliography}{1234}


\bibitem{Lo1} Agarwal S.,  Wettlaufer J.S.,
Maximal stochastic transport in the Lorenz equations.
Phys. Lett. A 380 (2016) 142-146.


\bibitem{fa} Applebaum D., L\'{e}vy Processes and Stochastic
Calculus. Second Edition. Cambridge University Press, 2009.


\bibitem{bd} Bao J., Mao X., Yin G., Yuan C.,
Competitive Lotka-Volterra population dynamics with jumps.
Nonlinear Anal. 74 (2011) 6601-6616.

\bibitem{be} Bao J., Yuan C., Stochastic population dynamics
driven by L\'{e}vy noise. J. Math. Anal. Appl. 391 (2012) 363-375.



\bibitem{a10} Chen F., Han Y., Li Y., Yang X., Periodic solutions
 of Fokker-Planck equations.
J. Differ. Equa. 263 (2017) 285-298.

\bibitem{CD} Chen L., Dong Z., Jiang J., Zhai J., On limiting behavior of stationary measures for stochastic evolution systems with small noise intensity. https://arxiv.org/pdf/1611.07223.pdf

 \bibitem{CCY} Chen X., Chen Z.Q., Tran K., Yin G.,  Recurrence and ergodicity for a class of regime-switching jump diffusions. Appl. Math. Optim. https://doi.org/10.1007/s00245-017-9470-9

\bibitem{3a} Da Prato G., Zabczyk J., Ergodicity for Infinite Dimensional
Systems. Cambridge University Press, 1996.

\bibitem{cc} Dong Y., Ergodicity of stochastic differential
equations driven by L\'{e}vy noise with local Lipschitz coefficients.
Adv. Math. 47 (2018) 11-47.

\bibitem{ck} Dong Z., On the uniqueness of invariant measure of the Burgers
equation driven by L\'{e}vy Processes. J. Theor. Probab. 21 (2008) 322-335.

\bibitem{DX} Dong Z., Xie Y., Ergodicity of stochastic 2D Navier-Stokes equation with
L\'evy noise. J. Differ. Equa. 251 (2011) 196-222.



\bibitem{a8} Khasminskii R.Z., Stochastic Stability of
Differential Equations. Springer-Verlag, Second Edition, 2012.

\bibitem{HL} Hu G., Li Y., Asymptotic behaviors of stochastic periodic differential equations with Markovian switching. Appl. Math. Compt. 264 (2015) 403-416.

\bibitem{a12} Hu H., Xu L., Existence and uniqueness theorems for periodic
Markov process and applications to stochastic functional differential
equations. J. Math. Anal. Appl. 466 (2018) 896-926.





\bibitem{a11} Li D., Xu D., Periodic solutions of stochastic
delay differential equations and applications to Logistic equation
and neural networks.
J. Korean Math. Soc. 50 (2013) 1165-1181.

\bibitem{lo} Lorenz E.N., Deterministic nonperiodic flow.
J. Atmos. Sci. 20 (1963) 130-141.

\bibitem{P125} Petersen K., Ergodic Theory. Cambridge University Press, 1983.





\bibitem{P} Poincar\'e H., Memoire sur les courbes definier par une equation differentiate. J. Math. Pures Appli. 3 (1881)
375-442; J. Math. Pures Appli. 3 (1882) 251-296; J. Math. Pures Appli. 4 (1885) 167-244; J. Math. Pures
Appli. 4 (1886) 151-217.

\bibitem{ma} Ren J., Wu J., Zhang X.,
Exponential ergodicity of multi-valued stochastic differential equations.
Bull. Sci. Math. France 134 (2010) 391-404.








\bibitem{ca} Xie B., Uniqueness of invariant measure of infinite
dimensional stochastic differential equations driven by L\'{e}vy noise. Potential Anal. 36 (2012) 35-66.

\bibitem{cj} Xie L., Zhang X., Ergodicity of stochastic differential equations with jumps and singular coefficients. https://arxiv.org/pdf/1705.07402.pdf

\bibitem{a9} Xu D., Huang Y., Yang Z., Existence theorems
for periodic Markov process and stochastic functional differential
equations.
Discrete Contin. Dyn. Syst. Ser. A 24 (2009) 1005-1023.



\bibitem{ra}Xu D., Li B., Long S., Teng L., Moment
estimate and existence for solutions of stochastic functional differential
equations.
Nonlinear Anal. 108 (2014) 128-143.

\bibitem{rb}Xu D., Li B., Long S., Teng L.,
Corrigendum to ``Moment estimate and existence for solutions of stochastic
functional differential equations".
Nonlinear Anal. 114 (2015) 40-41.

\bibitem{zg1} Zhang S., Meng X., Feng T., Zhang T.,
Dynamics analysis and numerical simulations of a stochastic
non-autonomous predator-prey system with
impulsive effects. Nonlinear Anal.: Hybrid Sys. 26 (2017) 19-37.

\bibitem{k1} Zhang X., Wang K., Li D.,
Stochastic periodic solutions of stochastic differential
equations driven by  L\'{e}vy process.
J. Math. Anal. Appl. 430 (2015) 231-242.


\end{thebibliography}
\end{document}